\documentclass[11pt]{article}

\usepackage{amsmath,amsfonts,amssymb}
\usepackage{graphicx}
\usepackage{color}
\usepackage{caption}
\usepackage{epsfig}
\usepackage{epstopdf}

\def\N{{\mathbb N}}
\def\R{{\mathbb R}}

\let\theta\vartheta
\let\phi\varphi

\let\union\cup

\def\vec#1{{\mathbf{#1}}}

\DeclareMathOperator{\sgn}{sgn}

\def\uh{\vec{u}_h}
\def\eh{\vec{e}_h}
\def\fh{\vec{f}_h}
\def\gh{\vec{g}_h}
\def\rh{\vec{r}_h}
\def\wh{\vec{w}_h}

\def\ehh{\vec{e}_{2h}}

\def\rhh{\vec{r}_{2h}}

\def\whh{\vec{w}_{2h}}

\title{A Finite Difference Ghost-cell Multigrid approach for Poisson Equation with mixed Boundary Conditions in Arbitrary Domain}
\author{Armando Coco, Giovanni Russo\thanks{Dipartimento di Matematica e Informatica, Universit\`a di Catania, Catania Italy}}

\begin{document}

\maketitle

\begin{abstract}
In this paper we present a multigrid approach to accelerate the convergence of the iterative method proposed in~\cite{CocoRusso:Elliptic} to solve the Poisson equation in arbitrary domain $\Omega$, identified by a level set function $\phi$, $\Omega=\left\{ x \in \mathbb{R}^d \colon \phi(x)<0 \right\}$, and mixed boundary conditions.
The discretization is based on finite difference scheme and ghost-cell method.
This multigrid strategy can be applied also to more general problems where a non-eliminated boundary condition approach is used. Arbitrary domain make the definition of the restriction operator for boundary conditions hard to find. A suitable restriction operator is provided in this work, together with a proper treatment of the boundary smoothing, in order to avoid degradation of the convergence factor of the multigrid due to boundary effects. Several numerical tests confirm the good convergence property of the new method.
\end{abstract}

%\tableofcontents
\section*{Introduction}
Multigrid technique is one of the most efficient strategy to solve a class of partial differential equations, using a hierarchy of discretizations. It accelerates the convergence of an existing iterative method, which otherwise slowly converges toward the solution of the discrete problem, due to the bad convergence rate for the low frequency components of the error. The idea of multigrid method is to solve such low frequency component in a coarser grid. An introduction to multigrid can be found, for example, in~\cite{Briggs:MG}, while more advanced textbook on the subject are, for example, ~\cite{Trottemberg:MG, Hackbusch:MG}.
Most iterative schemes to solve Elliptic equations can be speeded up by a multigrid technique.

Elliptic equation in arbitrary domain (possibly with moving boundary) is central to many applications, such as diffusion phenomena, fluid dynamics, charge transport in semiconductors, crystal growth, electromagnetism and many others. 
The wide range of applications may require different kind of boundary conditions. Let us look for instance at the temperature distribution in a medium of arbitrary shape satisfying stationary heat equation: we may have Dirichlet (the temperature is fixed at the boundary), Neumann (heat flux is prescribed), or mixed boundary conditions (namely different boundary conditions on different part of the boundary). More general Robin boundary conditions may also be prescribed, as in Stefan-type problem, in which a combination of temperature and heat flux is prescribed at the boundary (e.g. see~\cite{Gibou:fourth_order, Caflish:IslandDynamics}). An application we have in mind is to fluid dynamics: the aim is to model the motion of an incompressible fluid contained in a tank of arbitrary shape. The problem is modeled by incompressible Navier-Stokes equations, which are solved by projection method of Chorin~\cite{Chorin:projection, Chorin:projection1997}. This leads to an elliptic equation for the pressure, obtained enforcing the incompressibility condition. This pressure equation requires Dirichlet condition on the free surface of the fluid and Neumann condition on the rigid walls. The pressure equation is the bottleneck of the whole method and therefore requires an efficient solver.

Several techniques have been developed to solve Elliptic equation on an arbitrary domain. Finite Element Methods use a mesh triangulation to capture the boundary, such as in~\cite{Quart:PDE, QuartSacco:PDE, Nochetto:StefanProblemFEM,Schmidt:dendritiesFEM}. However, in presence of moving boundary, a grid re-meshing is needed at each time step, which makes the method expensive.
Furthermore, for a complex geometry, generation of a good mesh is a non trivial task that may require a considerable amount of work~\cite{Lonher:ImmersedEmbedded}.
For this reason a Cartesian grid method is preferred together with a level-set approach to keep track of the boundary at each time step. Level-set methods have been introduced to implicitly define a domain and its boundary, in order to simple handle complex topological changes of moving boundary such as merging and breaking up. Several papers and books exist in the literature about level-set method:~\cite{SSO:level_set, Osher-Fedkiw:level_set, Russo-Smereka:reconstruction, Gibou:reinizialization, Sethian:level_set} are just some examples.

Since the boundary may be not aligned with the grid, a special treatment is needed. The simplest method makes use of the Shortley-Weller discretization~\cite{Shortley-Weller:discretization}, that discretizes the Laplacian operator with usual central difference away from the boundary and makes use of a non symmetric stencil in the interior points of the domain close to the boundary. While this discretization provides a simple second order method for Dirichlet conditions, it cannot be immediately applied in presence of Neumann conditions. In fact, Shortley-Weller discretization~\cite{Shortley-Weller:discretization} for Neumann conditions requires that the value of the numerical solution is suitably reconstructed at the intersection between the grid and the boundary by applying the boundary condition. This approach is adopted, for example, by Hackbusch in~\cite{Hackbusch:elliptic} to first order accuracy, and by other authors (see~\cite{Bramble:ell} and the references therein) to second order accuracy. However, the method proposed by Bramble in~\cite{Bramble:ell} for second order accuracy is quite involved and not recommendable for practical purposes.

Another class of methods is based on cut-cell methods, obtained by a Finite Volume discretization which embeds the domain in a regular Cartesian grid~\cite{Colella:PoissonFV}. Cells that are cut by the boundary requires a special treatment, such as cell-merging and rotated-cell, in order to avoid a too strict restriction of the time step dictated by the CFL condition (e.g. see~\cite{Berger:rotated, vanLeer:CFLCutCell, Clarke:CFLCutCell}).

Other methods for Dirichlet condition are the Immersed Boundary Method, first proposed by Peskin in~\cite{Peskin:immersedInterface}, and later developed by several other authors~\cite{LeVequeLi:IIM, LeVequeLiWiegmann:CRACK}, with a proper multigrid approach~\cite{Li:IIM_MG}, and penalization methods~\cite{Iollo:penalization}.

In our method we will use a rather simple finite-difference ghost-cell technique, that adds extra grid points (ghost points) outside the domain in order to keep unchanged the symmetry of the stencil even for inside points close to the boundary.
A detailed description of the method can be found in~\cite{CocoRusso:Elliptic}.

In ghost points the boundary conditions are enforced in order to close the discrete system. The ghost-cell method was first developed by Fedkiw in~\cite{Fedkiw:GFM}, where a two-phase contact discontinuity was discretized (Ghost Fluid Method).
A second-order accurate method for Dirichlet conditions on regular Cartesian grid is proposed by Gibou \emph{et al.} in~\cite{Gibou:Ghost}. The value at the ghost nodes is assigned by linear extrapolation, and the whole discretization leads to a symmetric linear system, easily solved by a preconditioned conjugate gradient method. A fourth order accurate method is also proposed in~\cite{Gibou:fourth_order}. Other methods use a non-regular Cartesian grid, such as in~\cite{Gibou:quadtree}, where Gibou \emph{et al.} present finite difference schemes for solving the variable coefficient Poisson equation and heat equation on irregular domains with Dirichlet boundary conditions, using adaptive Cartesian grids. One efficient discretization based on cut-cell method to solve more general Robin conditions is proposed by Gibou \textit{et al.} in~\cite{Gibou:Robin}, which provides second order accuracy for the Poisson and heat equation and first order accuracy for Stefan-type problems.

Most of the techniques listed above cannot be straightforwardly applied in the special case of mixed boundary conditions. For cut-cell based methods~\cite{Colella:PoissonFV, Gibou:Robin}, different boundary conditions cannot to be imposed on the same boundary edge of a cut cell.
Simple efficient methods based on symmetric image of ghost points to solve mixed boundary condition problems provided with a multigrid algorithm have been recently developed in~\cite{Catalano:MG} and by Ma \textit{et al.} in~\cite{Ma:MG}.  

In our method~\cite{CocoRusso:Elliptic}, boundary conditions are neither eliminated from the discrete system (they are strongly coupled and their elimination is too hard to perform in more than one dimension) nor directly enforced (which leads to a non-convergent iterative method): they are \emph{relaxed together with the interior equations}. This leads us to an iterative scheme for the set of all unknowns (internal points and ghost points), which is proved to converge, at least for first order accurate discretization.

In this paper we provide a general multigrid technique to solve the discrete system coming from a continuous elliptic problem in case of non-eliminated boundary conditions. The smoothing procedure of the multigrid approach in the interior is Gauss-Seidel-like, while the iterations on the boundary are performed in order to provide smooth errors.

Multigrid techniques for non-eliminated boundary conditions are well-studied in literature in the case of rectangular domain (as we can see in~\cite{Hackbusch:MG, Trottemberg:MG}), where a restriction operator is defined separately for the interior of the domain and for the boundary, and the restriction of the boundary is performed using a restriction operator of codimension $1$, since ghost points are aligned with the Cartesian axis. In the case of arbitrary domain, ghost points have an irregular structure and we provide a reasonable definition of the restriction operator for the boundary conditions.
%Multigrid techniques in arbitrary domain (without using the ghost method) are proposed in several papers, such as in~\cite{Li:IIM_MG, Wan:MG}, where the case of discontinuous coefficient across an interface is treated. A survey of multigrid method for discontinuous coefficient is presented in~\cite{Chan:surveyMG}.
The method proposed in this paper can be extended to the case of discontinuous coefficient: a preliminary result in one dimension can be found in~\cite{CocoRusso:Ischia2010}, while the two-dimensional case is in preparation.

In this paper we also show that a proper treatment of the boundary iterations can improve the rate of convergence of the multigrid, making it closer to the convergence rate predicted by the Local Fourier Analysis for inside equations, as suggested by Brandt in~\cite{Brandt:RigorousAnalysisMG}. The cost of this extra computational work is negligible, i.e. tends to zero as the dimension of the problem increases.
A comparison with other kinds of treatment of the boundary condition smoothing procedure is carried out.

The paper is divided in three sections. We start with the multigrid approach in the one dimensional case, described in Section 1, with a special treatment of the transfer operators. Most of this method can be extended to high dimension, treated in Section 2, but a special care has to be taken for transferring the defect of boundary conditions. In this section a level-set approach is also introduced. Section 3 provides a strategy to improve the convergence factor making it closer to the one predicted by the Local Fourier Analysis for interior relaxations. Numerical evidence of such improvement is provided, together with a comparison with other boundary condition smoothers (Kaczmarz and Block relaxation).

In all the paper, we mainly use the notation of~\cite{Trottemberg:MG}.

\section{One-dimensional case}\label{1Dsection}
In this section we will follow the description of the method proposed in~\cite{CocoRusso:Elliptic}, which is second order accurate, and provides a multigrid approach to speed up the convergence of the iterative scheme. For one-dimensional case, the multigrid approach in arbitrary domain is a natural extension of the basic multigrid strategy that can be found in any good basic text about multigrid, such as~\cite{Trottemberg:MG, Briggs:MG, Hackbusch:MG}. Although if we can eliminate the boundary conditions from the linear system obtained by discretizing the problem, we always want to treat the case of non-eliminated boundary conditions in order to straightforwardly extend the method to more than one dimension, where the elimination of the boundary conditions from the system is hard to perform.

\subsection{Model problem and relaxation scheme}
\begin{figure}[!hbt]%
  \begin{center}
  \includegraphics[width=10cm]{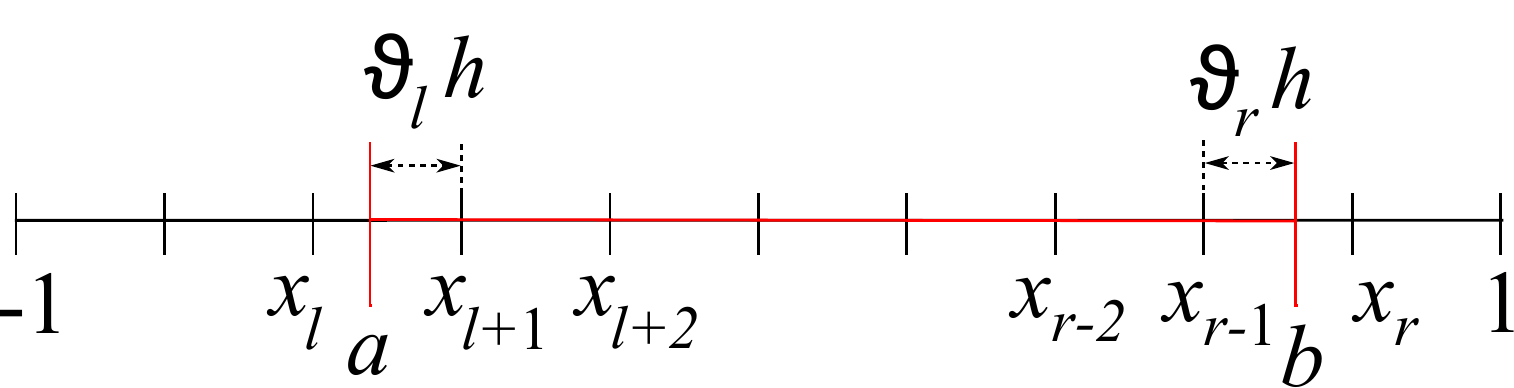}%
  \caption{\footnotesize{Discretization of the domain in 1D}}%
  \label{fig:Omega1d}%
  \end{center}
\end{figure}

Let $D=[-1,1]$ be the computational domain, $a$ and $b$ constants such that $-1<a<b<1$, and  $\Omega=[a,b]$. Let $N\geq 1$ be a fixed integer and $h=2/N$ the spatial step, let $D_h=\{-1=x_0<x_1< \ldots < x_N=1\}$ be the set of equally spaced grid points, and $\Omega_h = D_h \cap \Omega$ the set of inside grid points.
Consider the model problem
\begin{equation}\label{mainpr1d}
\begin{array}{rcll}
%\begin{cases}
-u'' &=&  f  & \text{ in } \Omega  \\
  u(a) &=& g_a & \\
   u'(b)   &=&  g_b. &
%\end{cases}
\end{array}
\end{equation}
Let $l$ and $r$ be such that $x_l \leq a < x_{l+1}$, $x_{r-1} < b \leq x_r$ 
(see Figure \ref{fig:Omega1d}). We use a ghost-cell method to discretize the problem. In order to obtain an iterative method, we solve the associate \textit{time-dependent} problem
\begin{eqnarray}
\frac{\partial u}{\partial t} & = & \frac{\partial^2 u}{\partial x^2} + f \quad \text{ in } \Omega  \label{mnprt1d}\\
\frac{\partial u(t, a)}{\partial t} & = & \mu_D \left(  g_a - u(t, a) \right)   \label{mnprd1d} \\
\frac{\partial u(t, b)}{\partial t} & = & \mu_N \left( g_b -  \frac{\partial u(t, b)}{\partial x} \right) \label{mnprn1d}\\
 u(0,x) & = & u_0(x)  \quad \text{ in } \Omega \label{mnpri1d}
\end{eqnarray}
and we look for a second order accurate steady state solution, which is the solution of the original problem.

Let us begin to discretize (\ref{mnprt1d}) in $\Omega_h \equiv \left\{ x_{l+1}, \ldots, x_{r-1} \right\}$.
We use central difference in space and forward Euler in time for (\ref{mnprt1d}) obtaining:
\begin{equation}\label{jact}
u_i^{(m+1)} = u_i^{(m)} + \frac{\Delta t}{h^2} \left( u_{i-1}^{(m)}-2 u_{i}^{(m)} +u_{i+1}^{(m)} \right) + \Delta t f_i, \; \; \; \; i= l+1,  \ldots, r-1.
\end{equation}
Taking the maximum time step consented by CFL condition~\cite{CFL:CFL}, i.e. $\Delta t=h^2/2$, we obtain:
\begin{equation}\label{jac}
u_i^{(m+1)} = \frac{1}{2} \left( u_{i-1}^{(m)}+u_{i+1}^{(m)} + h^2 f_i \right), \; \; \; \; i= l+1,  \ldots, r-1.
\end{equation}
Note that if we discretize directly the first equation (\ref{mainpr1d}) using central difference for the Laplacian operator, and use Jacobi iterative scheme for such discretization, we obtain exactly (\ref{jac}).

To obtain second order accuracy, we have to discretize the spatial terms in (\ref{mnprd1d}), (\ref{mnprn1d}) to second order, while first order time discretization can be used, because we are just interested at the accuracy as $t \rightarrow +\infty$.
%Let $\wh \colon \rightarrow \R$. We denote by
%\[
%\mathcal{L}_{x_j}^{(m)}+ \left(  \wh \right) \mbox{  and  } \mathcal{L}_{x_j}^{(m)}- \left(  \wh \right)
%\]
%the Lagrange interpolant of $\wh$ in nodes $x_j, \ldots, x_{j+n}$ and $x_j, \ldots, x_{j-n}$ respectively.

We then can use linear interpolation for $u(t,a)$ in (\ref{mnprd1d}). Since in some application it is required second order accuracy of the gradient of the solution, we use quadratic interpolation instead linear interpolation, obtaining: 
\begin{equation}\label{left2}
u_l^{(m+1)} = u_l^{(m)} - \mu_D \Delta t \left( (1+\theta_l) \frac{\theta_l}{2} u_l^{(m)} + (1+\theta_l)(1-\theta_l) u_{l+1}^{(m)} - (1-\theta_l) \frac{\theta_l}{2} u_{l+2}^{(m)} - g_a \right),
\end{equation}
where $\theta_l = (x_{l+1}-a)/h$, and quadratic interpolation of $u$ in nodes $x_r, x_{r-1}, x_{r-2}$ for the $\partial u(t, b) / \partial x$ in (\ref{mnprn1d}), obtaining:
\begin{eqnarray}
u_r^{(m+1)} &=& u_r^{(m)} - \frac{\mu_N \, \Delta t}{h} \left( u_{r-1}^{(m)} - u_{r-2}^{(m)} + \left(u_{r-2}^{(m)} - 2u_{r-1}^{(m)} + u_{r}^{(m)} \right)\left(\frac{1}{2} + \theta_r\right) \right) + \mu_N \, \Delta t \, g_b \label{right2}.
\end{eqnarray}

The constants $\mu_D$ and $\mu_N$ are chosen in order to satisfy the CFL conditions, i.e.\ $\mu_D \Delta t < 1$ and $\mu_N \Delta t /h < 2/3$ (see~\cite{CocoRusso:Elliptic}).
In numerical tests of Section \ref{numtests}, we choose $\mu_D \Delta t = 0.9$ and $\mu_N \Delta t = 0.9 \cdot 2h/(3)$.
Since $\Delta t = h^2/2$, then $\mu_D = 1.8/h^2$ and $\mu_N =  3.6/(3h)$.

In summary, our second order accurate iterative method is described by Eqs. 
(\ref{jac}), (\ref{left2}) and (\ref{right2}), with the choice of constants
\begin{equation}
\Delta t = h^2/2, \quad \mu_D = 1.8/h^2, \quad \mu_N = 1.2/h.
\label{eq:constants}
\end{equation}

%\textbf{Remark.} In some application it is required second order accuracy of the gradient of the solution instead just of the solution itself (see~\cite{CocoRusso:Elliptic}). In this case one should use quadratic interpolation in nodes $x_l$, $x_{l+1}$, $x_{l+2}$ for discretizing $u(t,a)$ in (\ref{mnprd1d}) in place of linear interpolation. 

\subsection{Multigrid approach}
We call $\Gamma_h$ the set of ghost points, i.e. $\Gamma_h = \{ x_l, x_r \}$. Let $I_h$ be a general subset of $D_h$. We introduce the linear space of grid functions over $I_h$ and we denote it $S(I_h) = \{ \wh \colon I_h \rightarrow \R \}$. For any $\wh \in S(I_h)$, we pose $w_i^h=\wh(x_i)$. Let $\fh \in S(\Omega_h)$ such that $f_i^h = f(x_i)$.
The iterative scheme (\ref{jac}), (\ref{left2}), (\ref{right2}) converges to the exact solution of the discretized system
\begin{eqnarray}
-\Delta_h \uh = \fh \label{mnpr1dh} \\
g^h_D(\uh) = g_a \label{mnprd1dh} \\
g^h_N(\uh) = g_b, \label{mnprn1dh}
\end{eqnarray}
where $\Delta_h \colon S(\Omega_h \cup \Gamma_h) \rightarrow S(\Omega_h)$ is defined by:
\[
\Delta_h \uh (x_i) = \frac{u_{i-1}^h - 2 u_i^h +u_{i+1}^h}{h^2}, \; \; \; x_i \in \Omega_h,
\]
while $g^h_D,g^h_N \colon S(\Omega_h \cup \Gamma_h) \rightarrow \R$ are the discrete versions of the boundary conditions:
\[
g^h_D(\uh) = (1+\theta_l) \frac{\theta_l}{2} u_l^{h} + (1+\theta_l)(1-\theta_l) u_{l+1}^{h} - (1-\theta_l) \frac{\theta_l}{2} u_{l+2}^{h},
\]
\[
g^h_N(\uh) = \frac{u_{r-1}^{h} - u_{r-2}^{h}}{h} + \frac{u_{r-2}^{h} - 2u_{r-1}^{h} + u_{r}^{h} }{h} \left(\frac{1}{2} + \theta_r\right).
\]
System (\ref{mnpr1dh})-(\ref{mnprn1dh}) can be interpreted in general as a discrete system of a Poisson equation with non-eliminated boundary conditions.

Let us consider an arbitrary grid function $\vec{\tilde{u}}^h \in S(\Omega_h \cup \Gamma_h)$ and let
\begin{eqnarray*}
\rh &=& \fh + \Delta_h \vec{\tilde{u}}^h \\
\tilde{g}_a &=& g_a - g^h_D(\vec{\tilde{u}}^h) \\
\tilde{g}_b &=& g_b - g^h_N(\vec{\tilde{u}}^h) 
\end{eqnarray*}
be the defects of (\ref{mnpr1dh}), (\ref{mnprd1dh}), (\ref{mnprn1dh}) respectively. Because of the linearity of $\Delta_h$, $g^h_D$, $g^h_N$, if we solve exactly the so-called \textit{residual problem}
\begin{eqnarray}
-\Delta_h \eh = \rh \\
g^h_D(\eh) = \tilde{g}_a \\
g^h_N(\eh) = \tilde{g}_b
\end{eqnarray}
in the unknown $\eh \in S(\Omega_h \cup \Gamma_h)$, then $\uh = \vec{\tilde{u}}^h + \eh$ is the exact solution of the system (\ref{mnpr1dh}), (\ref{mnprd1dh}), (\ref{mnprn1dh}). In the basic idea of multigrid one needs to solve the residual problem in a grid coarser than the original one.
\\
We can summarize the iterative scheme (\ref{jac}), (\ref{left2}), (\ref{right2}) as follows:
\begin{equation}\label{relax1d}
\uh^{(m+1)} = \Re_h \left( \uh^{(m)} , \fh, g_a, g_b \right)
\end{equation}
\begin{equation}\label{relaxO1d}
\Re_h \colon S(\Omega_h \cup \Gamma_h) \times S(\Omega_h) \times \R^2 \longrightarrow S(\Omega_h \cup \Gamma_h).
\end{equation}
Note that the iterative scheme  (\ref{jac}), (\ref{left2}), (\ref{right2}) is of a Jacobi kind. In order to provide a multigrid strategy, we just require that the iteration operator (\ref{relaxO1d}) has the \textit{smoothing property}, i.e. after few iteration steps (\ref{relax1d}), the error becomes smooth (not necessarily small). Roughly speaking, the high-frequency components of the error reduce quickly. We call \textit{smoothers} any operator (\ref{relaxO1d}) with this property.
Many iterators have this property, such as Gauss-Seidel or weighted Jacobi (with weight $\omega = 2/3$ in 1D or $\omega = 4/5$ in 2D), but not Jacobi (see~\cite[pag. 30--32]{Trottemberg:MG} for more details).
From now on, by (\ref{relax1d}) we shall intend the Gauss-Seidel version of (\ref{jac}), (\ref{left2}), (\ref{right2}), i.e.:
\begin{eqnarray*}
u_l^{(m+1)} &=& u_l^{(m)} - \mu_D \Delta t \left( \theta_l u_l^{(m)} + (1-\theta_l) u_{l+1}^{(m)} - g_a \right) \\
u_i^{(m+1)} &=& \frac{1}{2} \left( u_{i-1}^{(m+1)}+u_{i+1}^{(m)} + h^2 f_i \right), \; \; \; \; i= l+1,  \ldots, r-1 \\
u_r^{(m+1)} &=& u_r^{(m)} + \mu_N \: \Delta t \: g_b \\
& & - \frac{\mu_N \Delta t}{h} \left( u_{r-1}^{(m+1)} - u_{r-2}^{(m+1)} + \left(u_{r-2}^{(m+1)} - 2u_{r-1}^{(m+1)} + u_{r}^{(m+1)} \right)\left(\frac{1}{2} + \theta_r\right) \right).
\end{eqnarray*}
In order to explain the multigrid approach, we just describe the two-grid correction scheme (TGCS), because all the other schemes, such as $V$-cycle, $W$-cycle, $F$-cycle or Full multigrid cycle, can be easily derived from it (see~\cite[Sections 2.4, 2.6]{Trottemberg:MG} for more details). The TGCS consists into the following algorithm:
\begin{enumerate}\label{alg1d}
\item Set initial guess $\uh = 0$ \\
\item Relax $\nu_1$ times on the finest grid:
for $k$ from $1$ to $\nu_1$ do
\[
\uh \colon \! \! = \Re_h \left( \uh , \fh, g_a, g_b \right)
\]
\item Compute the defects
\begin{eqnarray*}
\rh &=& \fh + \Delta_h \uh \\
\tilde{g}_a &=& g_a - g^h_D(\vec{\tilde{u}}^h) \\
\tilde{g}_b &=& g_b - g^h_N(\vec{\tilde{u}}^h)
\end{eqnarray*}
\item Transfer the defect $\rh$ to a coarser grid with spatial step $2h$ by a suitable \textit{restriction operator}
\[
\rhh = I_{2h}^{h} \left( \rh \right)
\]
\item Solve exactly the residual problem on the coarser grid
\begin{eqnarray}
-\Delta_{2h} \ehh = \rhh \\
g^{2h}_D(\ehh) = \tilde{g}_a \\
g^{2h}_N(\ehh) = \tilde{g}_b
\end{eqnarray}
in the unknown $\ehh \in S(\Omega_{2h} \cup \Gamma_{2h})$
\item Transfer the error to the finest grid by a suitable \textit{interpolation operator}
\[
\eh = I_h^{2h} \left( \ehh \right)
\] 
\item Correct the fine-grid approximation
\[
\uh \colon \! \!= \uh + \eh
\]
\item Relax $\nu_2$ times on the finest grid:
for $k$ from $1$ to $\nu_2$ do
\[
\uh \colon \! \! = \Re_h \left( \uh , \fh, g_a, g_b \right)
\]

%\item Compute the defect $\rh \in S(D_h)$:
%\[
%\rh(x,y) =
%\left\{
%\begin{array}{cc}
%\left( \fh + \Delta_h \uh \right)(x,y) & \text{ if } (x,y) \in \Omega_h \\
%\left( g-L_h \uh \right)(x,y) & \text{ if } (x,y) \in \Gamma_h \\
%0 & \text{ otherwise }
%\end{array}
%\right.
%\]
\end{enumerate}
We have just to explain the steps concerning grid migration (steps 4 and 6).

\subsection{Transfer grid operators}
In this section, we describe the transfer grid operators for vertex-centered grid. We observe that our approach is based on the discretization of the equations on the various grids (both for inner and ghost points). This approach is very different from algebraic multigrid. As a consequence, the interpolation and the restriction operators are not the transpose of each other.

\subsubsection{Restriction operator}
\begin{figure}[!hbt]
 \begin{minipage}[c]{0.99\textwidth}
   	\centering
   	%\captionsetup{width=0.70\textwidth}
		\includegraphics[width=0.60\textwidth]{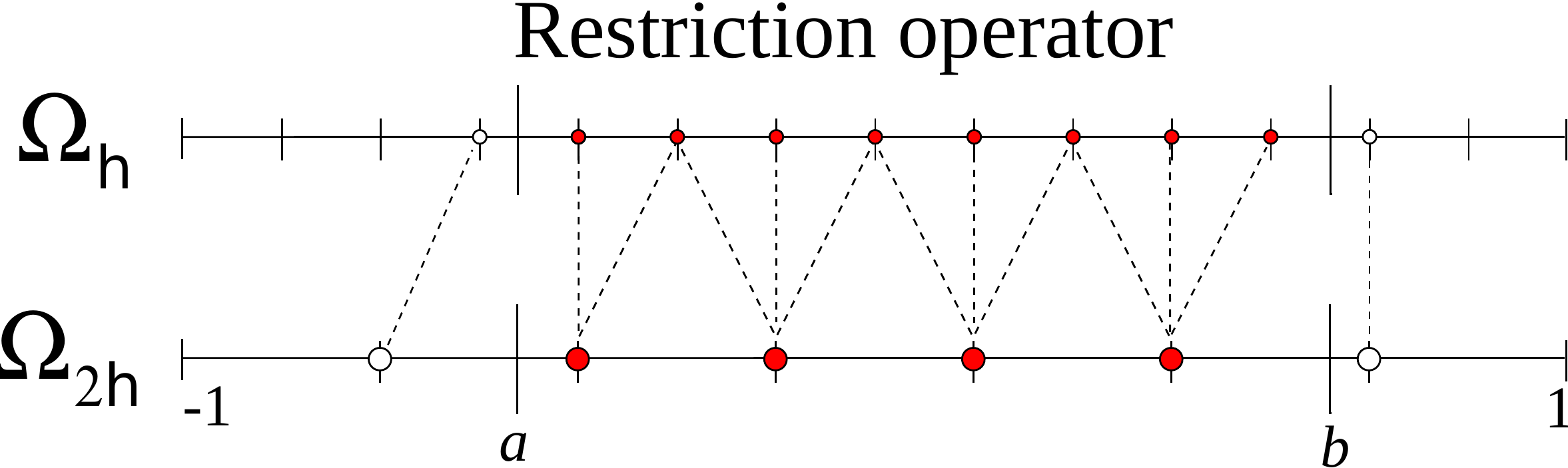}
	%\caption{\footnotesize{}}
	%\label{}
 \end{minipage}
 %\ \hspace{2mm} %\hspace{3mm} \
 \vskip 0.5 cm
 \begin{minipage}[c]{0.99\textwidth}
  	\centering
  	%\captionsetup{width=0.70\textwidth}
		\includegraphics[width=0.60\textwidth]{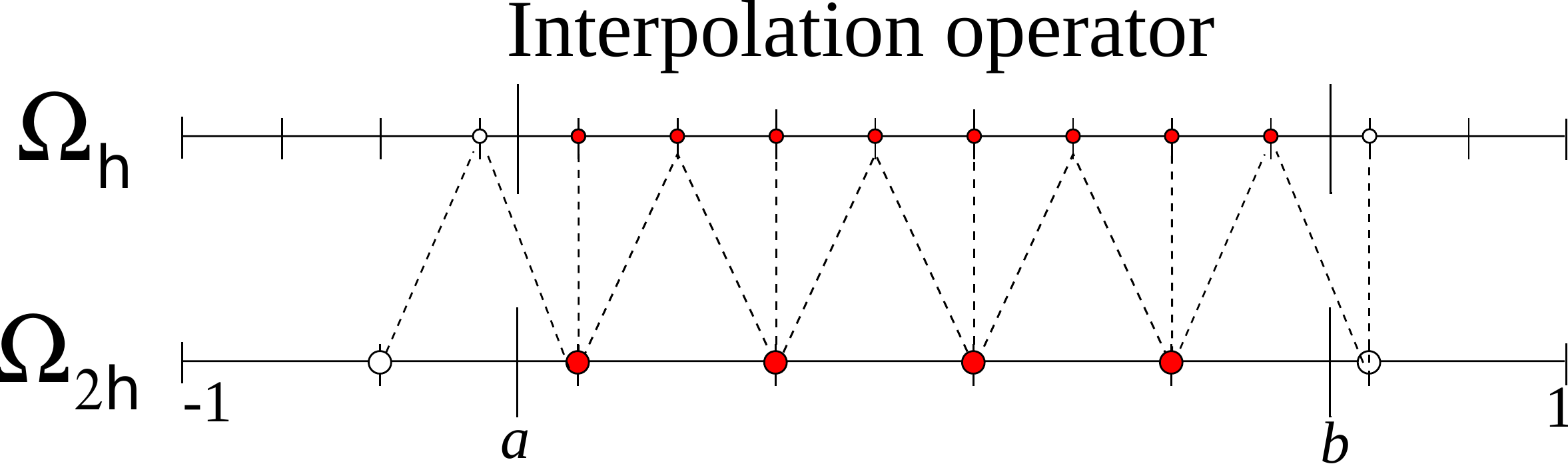}
	\caption{\footnotesize{ Vertex-centered discretization in 1D. Inner grid nodes (red circles) and ghost points (empty circles) on the fine and coarse mesh. The dashed lines represent the action of the restriction (up) and the interpolation (down) operators. } }
	\label{fig:VC}
 \end{minipage}
\end{figure}

%Non abbiamo descritto se usiamo vertex-centered o cell-centered. Dovremmo essere piu' specifici, e magari descrivere entrambe i metodi. Delle figure potrebbero accompagnare la descrizione.
Since such operator will act on the defect $\rh \in S(\Omega_h)$ (step 4), we must determine $I^h_{2h} \rh (x)$ for any $x \in \Omega_{2h}$ using only values inside $\Omega_h$. This is justified by the fact that the defect of the inside grid points (referred to the Poisson equation) may be very different (after few relaxations) from the defects $\tilde{g}_a$, $\tilde{g}_b$ (referred to the boundary conditions and stored computationally in the ghost points), because the operators (for inner equations and for boundary conditions) scale with different powers of $h$.
Then, let $x \in \Omega_{2h}$ and refer to Fig. \ref{fig:VC} (upper part). If $x$ is not near an outside grid point, i.e. $\min \{ \left| x- a \right|, \left| x- b \right| \} \geq h$, then we will use the standard full-weighting restriction operator (FW):
\begin{equation}\label{restI}
I^h_{2h} \rh (x) = \frac{1}{4} \left( \rh(x-h) + 2 \, \rh(x) + \rh(x+h) \right),
\end{equation}
while if $x-h < a$ or $x+h> b$ we set respectively
\begin{equation}\label{restL}
I^h_{2h} \rh (x) = \frac{1}{2} \left( \rh(x) + \rh(x+h) \right)
\end{equation}
or
\begin{equation}\label{restR}
I^h_{2h} \rh (x) = \frac{1}{2} \left( \rh(x-h) + \rh(x) \right).
\end{equation}

\subsubsection{Interpolation operator}
Since the interpolation operator acts on the error (step 6), which is continuous across the boundary, we do not need to separate the interpolation for inner equations from the interpolation of ghost points, and then we just use the standard linear interpolation operator (see the lower part of Fig. \ref{fig:VC}):
\[
\begin{array}{rcll}
I^{2h}_{h} \ehh (x_j) &=& \ehh(x_j) & \mbox{ if $j$ is even} \\
&&& \\
I^{2h}_{h} \ehh (x_j) &=& \frac{1}{2} \left( \ehh(x_{j-1}) + \ehh(x_{j+1}) \right) & \mbox{ if $j$ is odd.} 
\end{array}
\]

\textbf{Remark. 1 ($V$-cycle)} The $V$-cycle algorithm is easily obtained from the TGCS recursively, namely applying the same algorithm to solve the residual equation in step 5. To terminate the recursion, an exact solver is used to solve the residual problem when the grid becomes too coarse.

\textbf{Remark. 2 ($W$-cycle)} The $W$-cycle is similar to the $V$-cycle, with the only difference that the residual problem is solved recursively two times instead of one (in general schemes, $\gamma$ times, but $\gamma > 2$ is considered useless for practical purpose).

\textbf{Remark. 3 (Coarser operator)} We observe that the discrete operator $\Delta_{2h}$ in step 5 is just the operator obtained discretizing directly the continuous operator in the coarser grid, and not the operator obtained by the Galerkin condition
\[
\Delta_{2h} = I^{h}_{2h} \: \Delta_{h} \: I^{2h}_{h}.
\]
The latter approach, typical of algebraic multigrid, makes the algebraic problem more expensive from a computational point of view and does not take advantage of the fact that the discrete problem comes from a continuous problem.

\section{High-dimensional case}
%Io inizierei proprio dalla descrizione di un dominio rettangolare, i cui vertici on appartengono ad $\Omega_h$. In questo modo e' facile spiegare che la tecnica in 1D si trasferisce tale e quale, utilizzando il prodotto cartesiano. Poi passerei a domini piu' generali
In this case the defect of the boundary conditions has to be transferred in a suitable way to a coarser grid. The restriction has to be performed separately from the restriction of the interior equations, since these defects may show a sharp gradient crossing the boundary, because the discrete operators scale with different powers of $h$.

In case of arbitrary domain, ghost points may have a complex structure and the restriction cannot be defined straightforwardly as in the rectangular case, where ghost points are aligned with the grid and the restriction can be performed by a one dimensional operator.

For arbitrary domain we first need to extend the defect in a narrow band outside the domain constant along normal direction, and then we can operate the restriction as in the interior of the domain.
%In this case, the multigrid approach for rectangular domain cannot be easily extended for arbitrary domain. In order to bring back the problem as in a rectangular case, we will extend the defect outside the domain, in such a way it is defined in all the box in which the domain is embedded.
For the sake of clarity, we describe the multigrid strategy in the two-dimensional case, but the procedure can be extended straightforwardly in more dimensions.
We always refer to the second order method proposed in~\cite{CocoRusso:Elliptic}, which is briefly recalled here.

\subsection{Model problem and relaxation scheme}
Let $D = [-1,1]^2$ be the computational domain, $\Omega \subset D$ be a domain such that $\partial \Omega \cap \partial D = \emptyset$. Let ${\Gamma_D, \Gamma_N}$ be a partition of $\partial \Omega$ (i.e.\ $\Gamma_D \union \Gamma_N = \partial \Omega, \stackrel{\circ}{\Gamma}_D \cap \stackrel{\circ}{\Gamma}_N = \emptyset$, where the interior points are computed in the $d-1$ dimensional topological space). Consider the model problem
\begin{equation}\label{mainpr}
\begin{array}{rcll}
%\begin{cases}
- \Delta u &=& f  & \text{ in } \Omega \\
u &=& g_D & \text{ on } \Gamma_D  \\
\frac{\partial u}{\partial n} &=& g_N & \text{ on } \Gamma_N,
%\end{cases}
\end{array}
\end{equation}
where $\vec{n}$ is the outward unit normal, $\Delta = \displaystyle \frac{\partial^2 }{\partial x^2} + \displaystyle \frac{\partial^2 }{\partial y^2}$ is the Laplacian operator, $f \colon \Omega \rightarrow \R$, $g_D \colon \Gamma_D \rightarrow \R$, $g_N \colon \Gamma_N \rightarrow \R$ are assigned functions.

In order to solve the elliptic problem (\ref{mainpr}), we can transform it in an evolutive problem (with a fictitious time) that we call the associate \textit{time-dependent} problem:
\begin{eqnarray}
\frac{\partial \tilde{u}}{\partial t} = \Delta \tilde{u} + f   & \text{ in } \Omega \label{mnprt2d}\\
\frac{\partial \tilde{u}}{\partial t} =  \mu_D(g_D - \tilde{u}) & \text{ on } \Gamma_D \label{mnprd2d} \\
\frac{\partial \tilde{u}}{\partial t} = \mu_N\left(g_N - \frac{\partial \tilde{u}}{\partial n}\right) & \text{ on } \Gamma_N \label{mnprn2d}\\
\tilde{u} = \tilde{u}_0  & \text{ in } \Omega \text{, when } t=0 \label{mnpri2d},
\end{eqnarray}
where $\mu_D$ and $\mu_N$ are two positive constants. Then we look for the steady state solution. An iterative scheme can therefore be obtained by discretizing the associate \textit{time-dependent} problem and considering the time just as an iterative parameter.

\subsection{Level-set function}\label{levelset}
In order to keep track of the boundary $\Gamma$, we introduce the level set function $\phi_0 \colon D \rightarrow \R$, in such a way:
\[
(x,y) \in \stackrel{\circ}{\Omega} \Longleftrightarrow \phi_0(x,y) < 0, \; \;
(x,y) \in \partial \Omega \Longleftrightarrow \phi_0(x,y) = 0.
\]
The outward unit normal to the boundary is 
\[
\vec{n}=\frac{ \nabla \phi_0 }{ \left| \nabla \phi_0 \right| }.
\]
General references on the level set method for tracking interfaces are, for examples,~\cite{Osher-Fedkiw:level_set} or~\cite{Sethian:level_set}.
From the level set function, we can obtain the signed distance function $\phi$ by fast marching methods~\cite{Sethian:level_set} or by the reinitialization procedure based on the numerical solution of the following PDE 
\begin{equation}\label{sdf}
\frac{\partial \phi}{\partial t} = \sgn(\phi_0) \left(1- \left| \nabla \phi \right| \right),
\end{equation}
as we can see, for instance, in~\cite{SSO:level_set, Russo-Smereka:reconstruction, Gibou:reinizialization}.
A signed distance function is preferred to a simple level-set function because sharp gradients are avoided and it is simpler to compute the boundary closest point to a given ghost point.
Now we assume that $\left| \nabla \phi \right| = 1$ and suppose we know the signed distance function just at the grid nodes.
In practice, Eq.\ (\ref{sdf}) has to be solved for a few time steps, in order to compute the distance function a few grid points away from the boundary.
z
\subsection{Relaxation operator}
%\subsubsection*{Notation.}
Let us introduce some notation. Let $d \in \N$ be the dimension of the problem, $N\geq 1$ be an integer and $h=2/N$ the spatial step. Let $D_h = \vec{j}h, 
\vec{j} = (j_1,\ldots,j_d) \in \left\{-N, N\right\}^d$ and $\Omega_h = \Omega \cap D_h$ be the discrete versions of $D$ 
and $\Omega$ respectively.
$D_h$ is the set of \textit{grid points}.
Two points $\vec{x'}$, $\vec{x''}$ in $D_h$ are called \emph{neighbor} if $\sum_{j=1}^d \left| x'_j - x''_j \right|=h$.
We call \textit{ghost point} any grid point that is both outside $\Omega$ and neighbor to a grid point inside $\Omega$. We call $\Gamma_h$ the set of all ghost points. Let $I_h$ be a general subset of $D_h$. We introduce the linear space of grid functions over $I_h$ and we denote it $S(I_h) = \{ \wh \colon I_h \rightarrow \R \}$.

From now on, we shall consider $d=2$, but the results are valid also for $d>2$.

Then, we write the basic iterative scheme (\textit{relaxation scheme}) discretizing the \textit{time-dependent} problem (\ref{mnprt2d})-(\ref{mnpri2d}).
For any grid point $(jh,ih)$ of $\Omega_h$, we write an equation obtained from the discretization of (\ref{mnprt2d}) in such point,
using forward Euler in time and central difference in space and taking the maximum time step consented by the CFL condition, i.e. $\Delta t = h^2/4$ (in general it is $\Delta t = h^2/(2d)$):
\begin{equation}\label{jac2d}
u_{i,j}^{(m+1)} = 1/4 \left( h^2 f_{i,j} + u_{i-1,j}^{(m)} + u_{i+1,j}^{(m)} + u_{i,j-1}^{(m)} + u_{i,j+1}^{(m)} \right).
\end{equation}
Eq. (\ref{jac2d}) is equivalent to discretize directly the first equation of (\ref{mainpr}) using central difference in space and applying Jacobi iteration scheme. 

Since we have used the standard 5-point stencil even for grid point close to the boundary, we have added new unknowns to the linear system (\textit{ghost points}).

To close the system of equations (\ref{jac2d}), we must write one equation for each ghost point. This can be done in three simple steps. Let $G \equiv (x_G,y_G)$ be a ghost point.

\begin{figure}[!hbt]%
  \begin{center}
  \includegraphics[width=6cm]{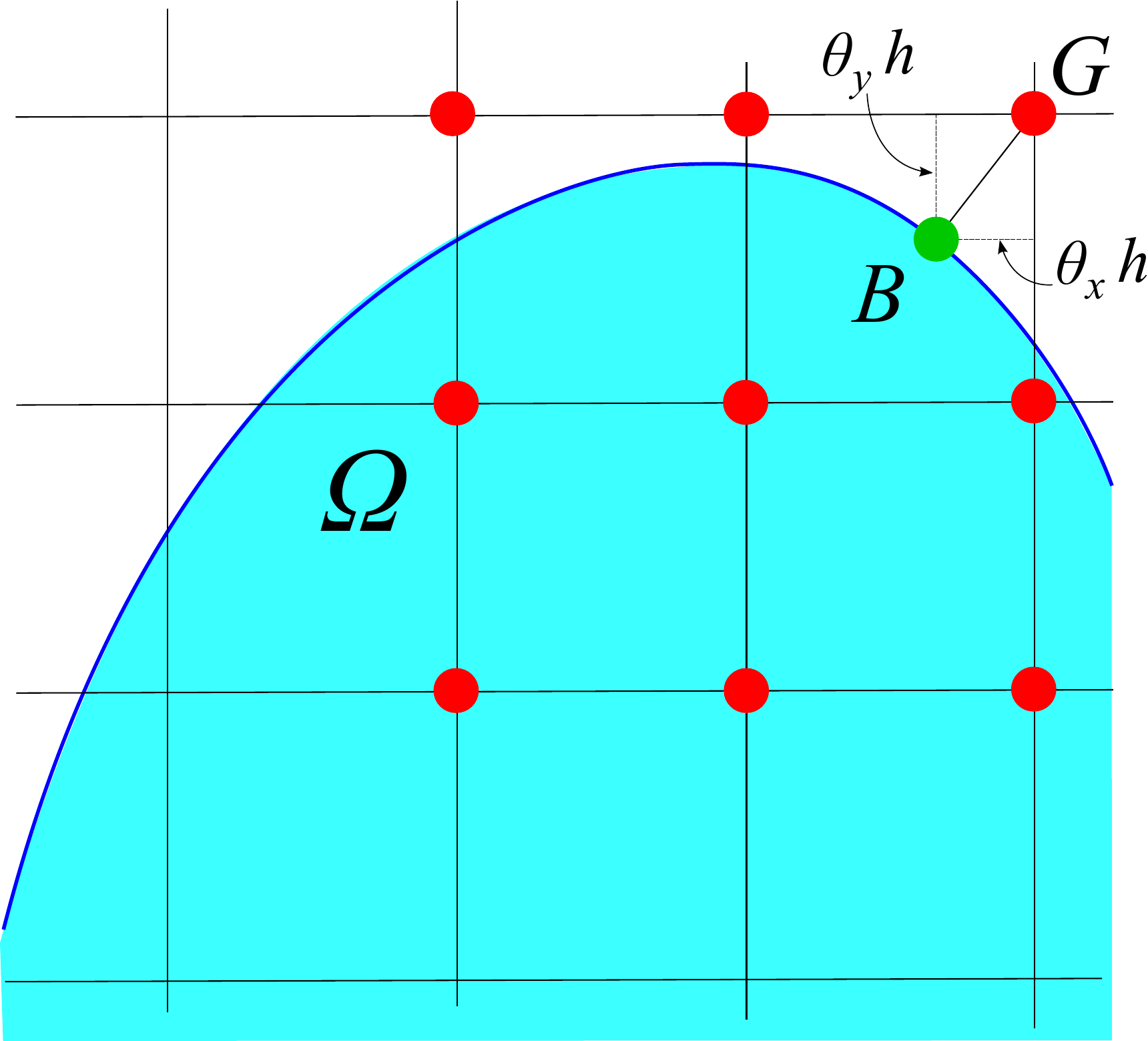}%
  \caption{\footnotesize{$B$ is the boundary closest point to $G$, while the red points are the nine-point stencil in Upwind direction referred to the ghost point $G$. }}%
  \label{fig:St9}%
  \end{center}
\end{figure}

\begin{enumerate}
\item Making use of the signed distance function $\phi$, we can compute the closest boundary point to $G$, that we call $B$ (see Figure \ref{fig:St9}), by:
\begin{equation}\label{clpt}
B \equiv (x_B,y_B) = G - \vec{n}_G \cdot \phi(G) = G - \left. \left( \frac{\nabla \phi }{ \left| \nabla \phi \right|} \right) \right|_G  \phi(G),
\end{equation}
using a second order accurate discretization for $\nabla \phi$, such as central difference in $G$.

\item Compute the nine-point stencil (depicted in Fig. \ref{fig:St9}) in Upwind direction, i.e.:
\begin{equation}\label{St_G}
St_G = \left\{ (x_G+ s_x \: k_1 \: h, y_G+ s_y \: k_2 \: h) \colon (k_1,k_2) \in \left\{ 0,1,2 \right\}^2 \right\},
\end{equation}
where $s_x=\sgn(x_B-x_G)$ and $s_y=\sgn(y_B-y_G)$.

\item Let $\mathcal{L}_{St_G}[u]$ be the biquadratic interpolant of the numerical solution $u$ in the stencil $St_G$.
If $B \in \Gamma_D$, the iteration for the ghost point $G$ will be obtained from the discretization of (\ref{mnprd2d}):
\begin{equation}\label{dir2d}
u_G^{(m+1)} = u_G^{(m)} + \mu_D\Delta t \left( g_D(B) - \mathcal{L}_{St_G} [u^{(m)}](B) \right)
\end{equation}
while if $B \in \Gamma_N$, the iteration for the ghost point $G$ will be obtained from the discretization of (\ref{mnprn2d}):
\begin{equation}\label{neu2d}
u_G^{(m+1)} = u_G^{(m)} + \mu_N\Delta t \left( g_N(B) -
\left. \left( \nabla \mathcal{L}_{St_G} [u^{(m)}] \cdot \frac{\nabla \mathcal{L}_{St_G} [\phi]}{\left| \nabla \mathcal{L}_{St_G} [\phi] \right|} \right) \right|_B \right).
\end{equation}
\end{enumerate}

The constants $\mu_D$ and $\mu_N$ are chosen in order to satisfy a CFL condition, i.e. $\mu_D \Delta t <1$ and
$\mu_N \Delta t < 2h/(3\sqrt{2})$~\cite{CocoRusso:Elliptic}.
In numerical tests of Section \ref{numtests}, we choose $\mu_D \Delta t = 0.9$ and $\mu_N \Delta t = 0.9 \cdot 2h/(3\sqrt{2})$.
Since $\Delta t = h^2/4$, then $\mu_D = 3.6/h^2$ and $\mu_N =  7.2/(3\sqrt{2}h)$.

\textbf{Remark 1 (Accuracy of \eqref{clpt}).} The accuracy of the evolution of point $B$ in (\ref{clpt}) depends on the accuracy at which $\phi$ is computed. If $\phi$ is known to order $h^p$, $p \in \left\{ 2,3 \right\}$, then $B$ will be computed to the same order of accuracy, provided we are far from singularities in $\phi$.

In Eq.\ (\ref{clpt}) one could omit the term $\left| \nabla \phi \right|$ in the denominator, because $\left| \nabla \phi \right|=1$ if $\phi$ is a signed distance function. However, it is better to keep such term, in case $\phi$ is only approximately a signed distance function.

\textbf{Remark 2 (Upwind stencil).} The reason for which we use an Upwind stencil $St_G$ is simple. Let us rewrite the Neumann boundary condition (\ref{mnprn2d}) as:
\begin{equation}\label{mnprn2dIP}
\frac{\partial u}{\partial t} + \mu_N \vec{n} \cdot  \nabla u = \mu_N  g_N,
\end{equation}
in such a way it appears to be an hyperbolic equation which propagates the solution $u$ along the characteristic (normal direction to the boundary) with speed $\mu_N \: \vec{n}$. Then, it is preferred to use an Upwind stencil in direction $-\vec{n}$ to discretize the spatial term, in order to guarantee convergence (see~\cite{Strikwerda:FD} for more detail about Upwind schemes in Conservation laws).

\textbf{Remark 3 (Reduced stencil).} If $St_G$ is not fully contained in $\Omega_h \cup \Gamma_h$, therefore we reduce the nine-point stencil to a smaller stencil, such as a $2 \times 2$ stencil or a (less accurate, more robust) three-point stencil. Such a reduction occurs rarely, and does not degrade the whole accuracy of the method (see~\cite{CocoRusso:Elliptic}).

Using the simplified notation, the iterative scheme converges to the solution of the problem:
%We can now summarize the discrete system described so far using the discrete operator notation, such as:
\begin{equation}\label{discr}
\begin{cases}
- \Delta_h \uh = \fh  \\
L_h \uh = \gh
\end{cases} 
\end{equation}
where:
\begin{itemize}
\item $\uh \in S(\Omega_h \cup \Gamma_h)$ is the unknown; 
\item
$\Delta_h \colon S(\Omega_h \cup \Gamma_h) \rightarrow S(\Omega_h)$ is the standard discrete version of the Laplacian  operator, namely:
\[
\Delta_h \wh (x,y) = \frac{1}{h^2} \left( \wh(x+h,y)+\wh(x-h,y)-4\wh(x,y)+\wh(x,y+h)+\wh(x,y-h) \right)
\]
for any $\wh \in S(\Omega_h \cup \Gamma_h)$ and $(x,y) \in \Omega_h$;
\item
$\fh \in S(\Omega_h)$ is defined by $f_h(P) = f(P)$ for any grid point $P \in \Omega_h$;
\item
$L_h \colon S(\Omega_h \cup \Gamma h) \rightarrow S(\Gamma_h)$ is the discrete version of boundary conditions, namely:
\begin{equation}\label{discrBC}
L_h \wh(G) =
\left\{
\begin{array}{ccc}
\mathcal{L}_{St_G}[u](B) & \mbox{ if } & B \in \Gamma_D \\
\left. \left( \nabla \mathcal{L}_{St_G} [u] \cdot \frac{\nabla \mathcal{L}_{St_G} [\phi]}{\left| \nabla \mathcal{L}_{St_G} [\phi] \right|} \right) \right|_B & \mbox{ if } & B \in \Gamma_N
\end{array}
\right.
\end{equation}
for any $\wh \in S(\Omega_h \cup \Gamma_h )$ and $G \in \Gamma_h$;
\item
$\gh \in S(\Gamma_h)$ is defined by:
\[
\gh(G) =
\left\{
\begin{array}{ccc}
g_D(B) & \mbox{ if } & B \in \Gamma_D \\
g_N(B) & \mbox{ if } & B \in \Gamma_N
\end{array}
\right.
\]
for any ghost point $G \in \Gamma_h$.
\end{itemize}

\subsection{Multigrid approach in 2D}
Consider the Poisson problem (\ref{mainpr}) and suppose we have a discrete approximation of the form (\ref{discr}). Therefore, we are dealing with non-eliminated boundary conditions. Let us introduce, for any spatial step $h$, an exact solver
\begin{equation}\label{solver}
\uh = S_h \left( \fh, \gh \right)
\end{equation}
of the system (\ref{discr}), and denote by
\begin{equation}\label{relaxO}
\Re_h \colon S(\Omega_h \cup \Gamma_h) \times S(\Omega_h) \times S(\Gamma_h) \longrightarrow S(\Omega_h \cup \Gamma_h)
\end{equation}
the relaxation operator, namely the iterative scheme
\begin{equation}\label{relax}
\uh^{(m+1)} = \Re_h \left( \uh^{(m)} , \fh, \gh \right)
\end{equation}
converges to the solution of (\ref{discr}) as $n \rightarrow + \infty$. In details, the iteration (\ref{relax}) summarize the iterative scheme (\ref{jac2d}), (\ref{dir2d}), (\ref{neu2d}).

As in one dimensional case, we will intend by (\ref{relax}) the Gauss-Seidel version of (\ref{jac2d}), (\ref{dir2d}), (\ref{neu2d}), in order to deal with a proper \textit{smoother},
and we have to order all points of $\Omega_h \cup \Gamma_h$ in some way. Let us choose the lexicographic ordering (GS-LEX):
\[
(x',y') \leq (x'',y'') \Longleftrightarrow x'<x'' \mbox{ or } x'=x'', \; y'<y''.
\]
The relaxation scheme can be easily extended to more efficient kinds of smoothers, such as Red-Black Gauss-Seidel (see~\cite{Trottemberg:MG}): however, we limit ourselves to study of the GS-LEX smoother.

In order to explain the multigrid approach, we just describe the two-grid correction scheme (TGCS), because all the others schemes, such as $V$-cycle, $W$-cycle or Full multigrid, can be easily derived from it (see~\cite[Sections 2.4, 2.6]{Trottemberg:MG} for more details).
The TGCS consists into the following algorithm:
\begin{enumerate}\label{alg}
\item Set initial guess $\uh = 0$ \\
\item Relax $\nu_1$ times on the finest grid:
for $k$ from $1$ to $\nu_1$ do
\[
\uh \colon \! \! = \Re_h \left( \uh , \fh, \gh \right)
\]
\item Compute the defect $\rh \in S(\Omega_h \cup \Gamma_h)$ such as:
\[
\rh = \left\{
 \begin{array}{ll}
  × f_h+\Delta_h \, \uh & \mbox{ in } \Omega_h \\
  × g_h - L_h \, \uh & \mbox{ on } \Gamma_h
 \end{array}
\right.
\]
\item Transfer the defect to a coarser grid with spatial step $2h$ by a suitable \textit{restriction operator}
\[
\rhh = I_{2h}^{h} \left( \rh \right)
\]
\item Solve exactly the residual problem in the coarser grid
\[
\ehh = S_{2h} \left( \rhh^I, \rhh^G \right)
\]
where $\rhh^I = \left. \rhh \right|_{\Omega_{2h}}$ and $\rhh^G = \left. \rhh \right|_{\Gamma_{2h}}$;
\item Transfer the error to the finest grid by a suitable \textit{interpolation operator}
\[
\eh = I_h^{2h} \left( \ehh \right)
\] 
\item Correct the fine-grid approximation
\[
\uh \colon \! \! = \uh + \eh
\]
\item Relax $\nu_2$ times on the finest grid:
for $k$ from $1$ to $\nu_2$ do
\[
\uh \colon \! \! = \Re_h \left( \uh , \fh, \gh \right)
\]
%\item Compute the defect $\rh \in S(D_h)$:
%\[
%\rh(x,y) =
%\left\{
%\begin{array}{cc}
%\left( \fh + \Delta_h \uh \right)(x,y) & \text{ if } (x,y) \in \Omega_h \\
%\left( g-L_h \uh \right)(x,y) & \text{ if } (x,y) \in \Gamma_h \\
%0 & \text{ otherwise }
%\end{array}
%\right.
%\]
\end{enumerate}
We have just to explain the steps concerning grid migration (steps 4 and 6). All the other steps are clear.
%defined in the previous section.
\subsection{Transfer grid operators}
If (\ref{relaxO}) has the smoothing property, after $\nu_1$ relaxations (step 2 of the algorithm) we have a smooth defect $\rh$. Therefore, we can hope to transfer this defect to a coarser grid without losing much information. The defect $\rh$ defined in the step 3 belongs to $S(\Omega_h \cup \Gamma_h)$. In order to transfer it to a coarser grid, it is convenient to extend in some way this defect in the whole computational domain $D_h$ (i.e. $\rh \in S(D_h)$), in such a way we can use the standard full-weighting stencil for the restriction operator $I_{2h}^h \colon S(D_h) \rightarrow S(D_{2h})$, that is (see~\cite[pag. 42]{Trottemberg:MG})
\begin{equation}\label{Fstencil}
I_{2h}^h = \frac{1}{16}
\left[
\begin{array}{ccc}
1 & 2 & 1 \\
2 & 4 & 2 \\
1 & 2 & 1 \\
\end{array}
\right]_{2h}^{h}.
\end{equation}
In general, by the stencil notation
\[
I_{2h}^h = 
\left[
\begin{array}{ccccc}
 & \vdots & \vdots & \vdots & \\
\cdots & t_{-1,-1} & t_{-1,0} & t_{-1,1} & \cdots \\
\cdots & t_{0,-1} & t_{0,0} & t_{0,1} & \cdots \\
\cdots & t_{1,-1} & t_{1,0} & t_{1,1} & \cdots \\
       & \vdots   &  \vdots &  \vdots &      
\end{array}
\right]_{2h}^{h}
\]
we will intend the restriction operator $I_{2h}^h$ defined by:
\[
   I_{2h}^h \wh (x,y) = \sum_{(i,j)\in R_k} t_{i,j} \wh(x+jh,y+ih),
\]
where only a finite number of coefficients $t_{i,j}$ is different from zero, and $R_k\equiv \left\{ -k,\ldots,k \right\}^2$ for some positive integer $k$. In practice $k=1$ allows second order restriction operator.

\begin{figure}[!hbt]%
  \begin{center}
  \includegraphics[width=0.65\textwidth]{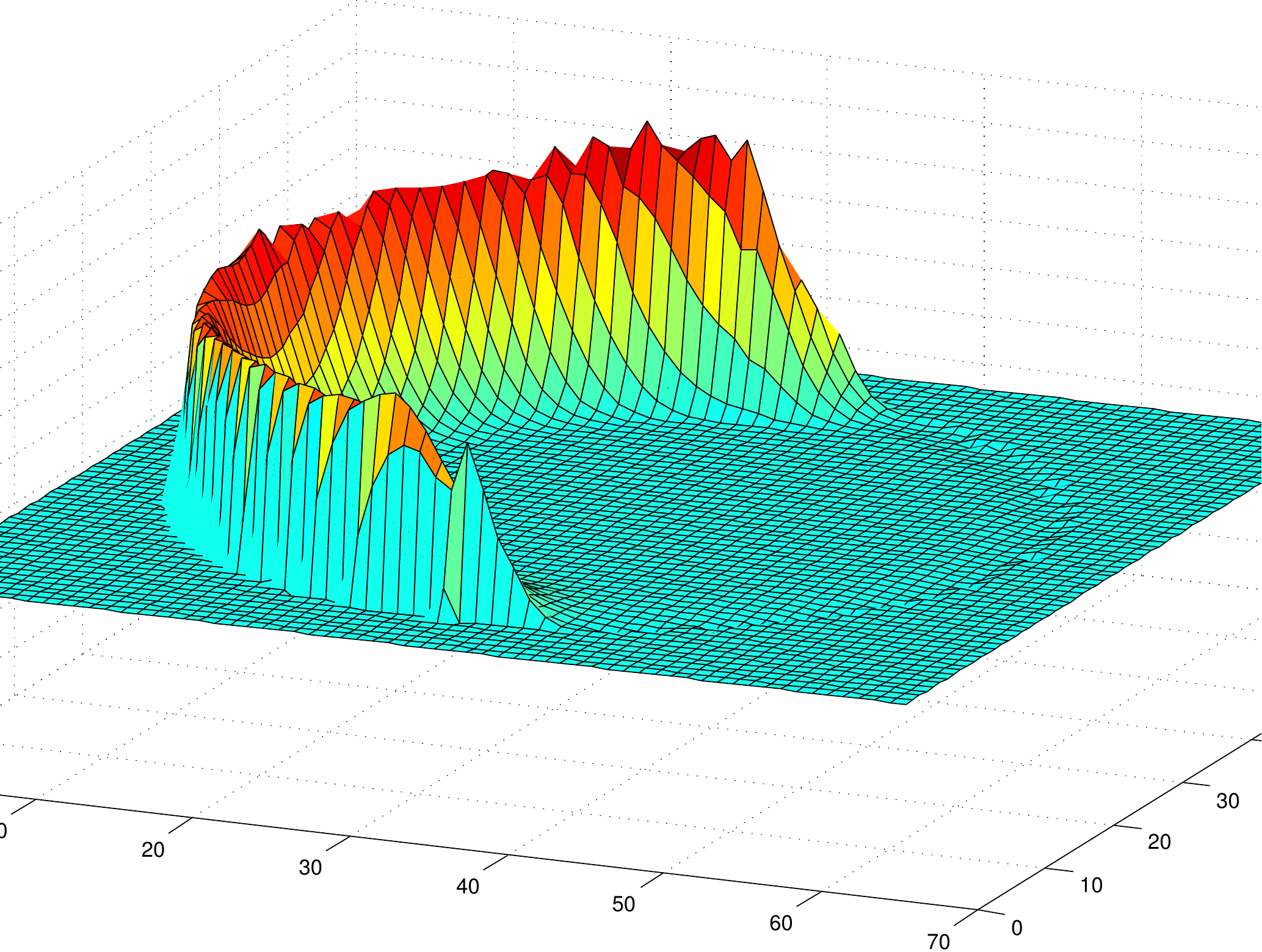}%
  \caption{\footnotesize{After few relaxations, the defect show a sharp gradient crossing the boundary (the defect is not yet extended outside), especially for the Neumann boundary. Here the domain is a circle.}}%
  \label{fig:resG}%
  \end{center}
\end{figure}

Let us suppose we have extended the defect to the whole computational domain $D_h$ (as it is carefully described in Sec.\ \ref{extension}). Anyhow, since we have different operators for inner equations and for boundary conditions, the defect is smooth separately inside $\Omega_h$ and along the ghost point $\Gamma_h$ (or $D_h - \Omega_h$, because of the extension), but it is not smooth in all $\Omega_h \cup \Gamma_h$ (it shows a sharp gradient crossing the boundary, as we can see in Fig. \ref{fig:resG}).

For this reason, it is convenient to transfer separately on the coarse grid the defect in $\Omega_h$ and in $D_h - \Omega_h$. To do that, we introduce a partial grid transfer
\[
\tilde{I}_{2h}^h \colon S(D_h) \times \mathcal{P}(D_h) \times S(D_{2h}) \longrightarrow S(D_{2h}),
\]
where $\mathcal{P}(D_h)$ is the family of all subset of $D_h$. Roughly speaking,
\[
\whh = \tilde{I}_{2h}^h \left( \wh, I_h, \tilde{w}_{2h} \right)
\]
means that we transfer $\wh$ to a coarser grid $\whh$ using only the points of $I_h$, leaving unaltered the value in the points out of $I_h$ and already stored in $\tilde{w}_{2h}$ (to better understand, we can think $I_h=\Omega_h$).

In details, let $(x,y) \in D_{2h} \cap I_h$. We focus our attention to the neighborhood of $(x,y)$, that is $\mathcal{N}(x,y) = \{ (x+jh,y+ih) \colon  j,i = -1,0,1 \}$.

\begin{figure}[!hbt]%
\begin{minipage}{0.30\textwidth}
  \begin{center}
  \includegraphics[width=1.00\textwidth]{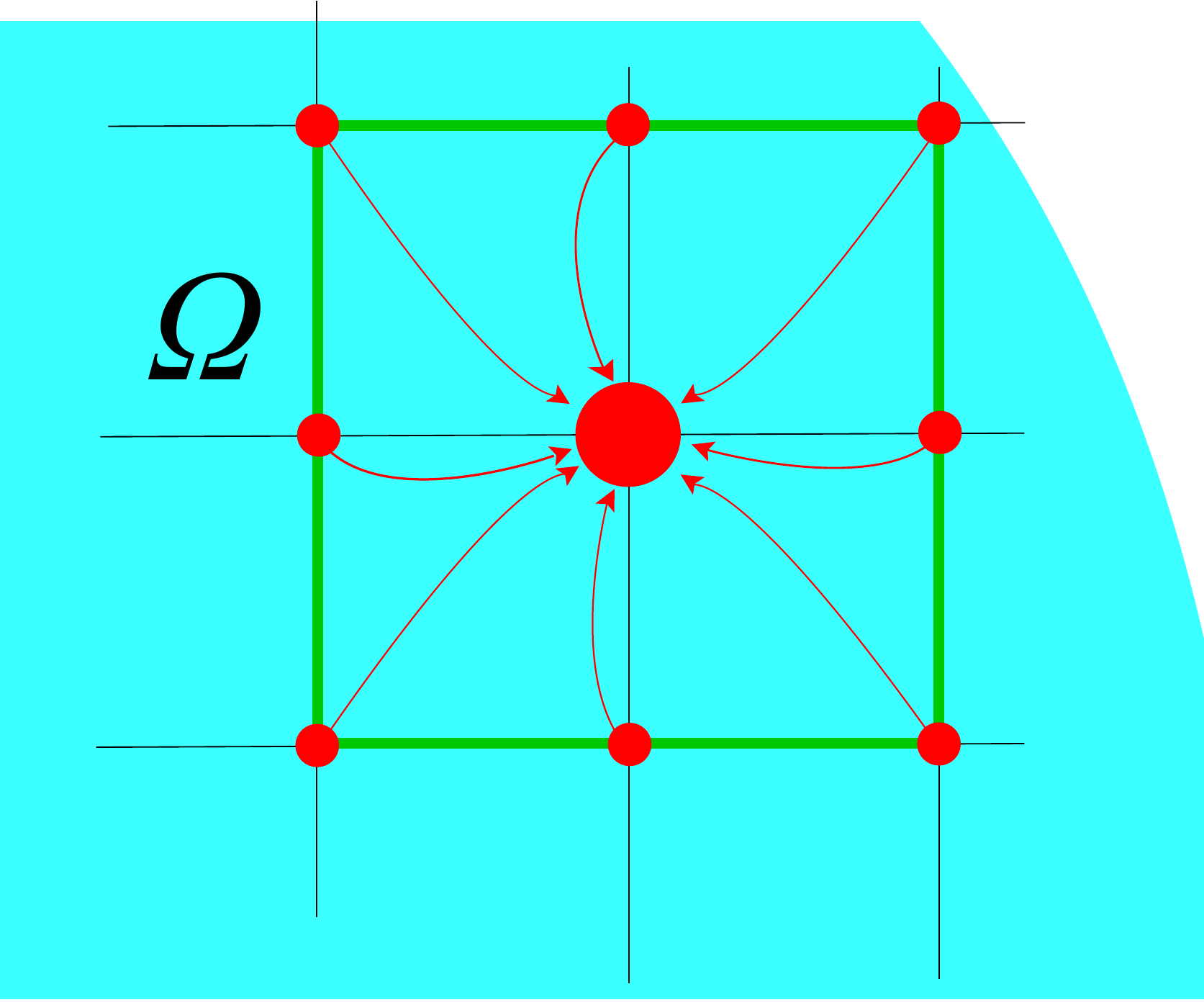}%
    %\caption{}%
  \end{center}
\end{minipage}
\hspace{1mm}
\begin{minipage}{0.30\textwidth}
  \begin{center}
  \includegraphics[width=1.00\textwidth]{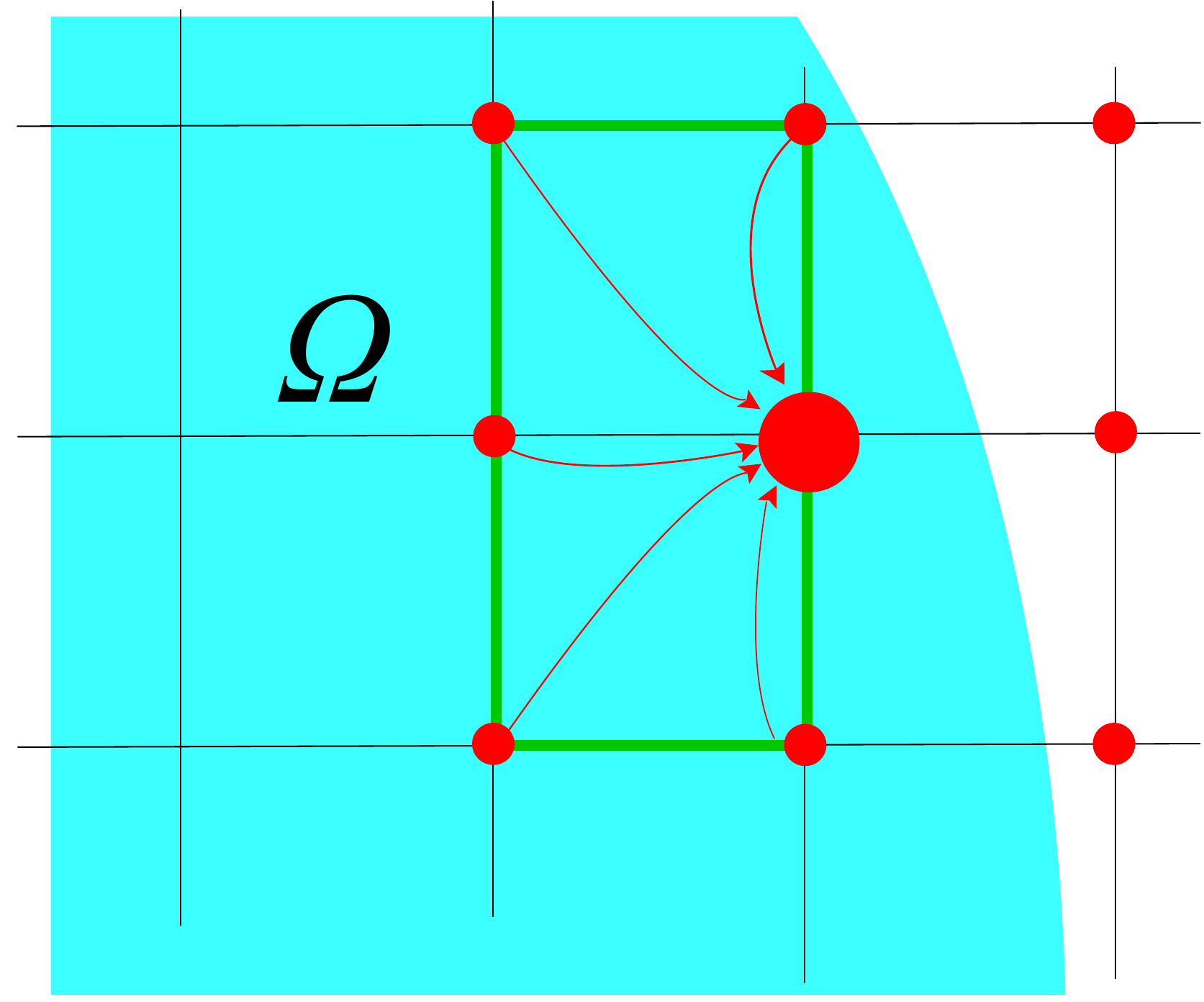}%
   %\caption{}
  \end{center}
\end{minipage}
\hspace{1mm}
\begin{minipage}{0.30\textwidth}
  \begin{center}
  \includegraphics[width=1.00\textwidth]{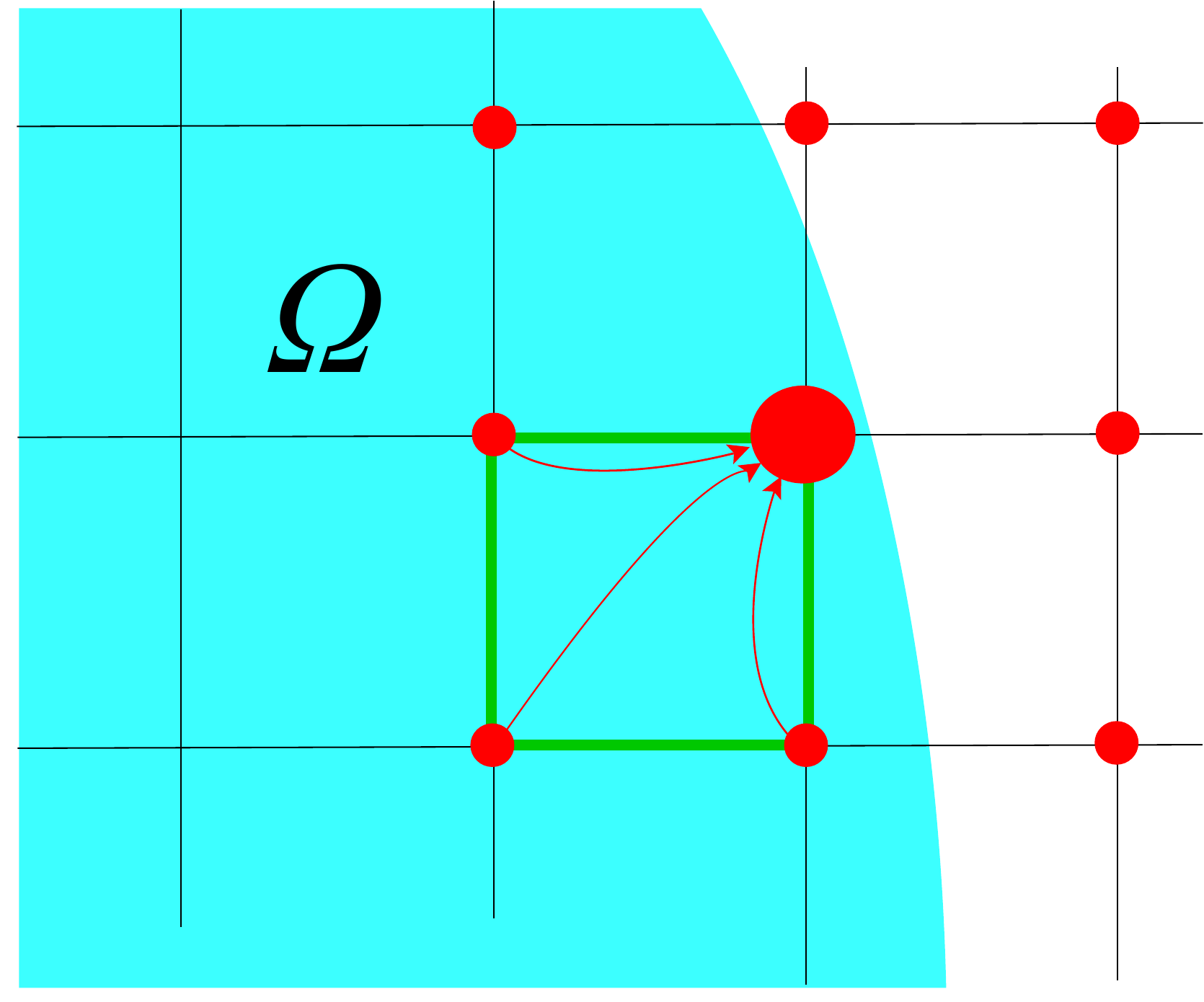}%
   %\caption{}
  \end{center}
\end{minipage}
 %\label{fig:Neigh}
%\vskip 0.1 cm
\begin{minipage}{0.30\textwidth}
\begin{equation*}
\frac{1}{16}
\left[
\begin{array}{ccc}
1 & 2 & 1 \\
2 & 4 & 2 \\
1 & 2 & 1 \\
\end{array}
\right]_{2h}^{h}
\end{equation*}
\end{minipage}
\hspace{1mm}
\begin{minipage}{0.30\textwidth}
\begin{equation*}
\frac{1}{16}
\left[
\begin{array}{ccc}
2 & 2 & 0 \\
4 & 4 & 0 \\
2 & 2 & 0 \\
\end{array}
\right]_{2h}^{h}
\end{equation*}
\end{minipage}
\hspace{1mm}
\begin{minipage}{0.30\textwidth}
\begin{equation*}
\frac{1}{16}
\left[
\begin{array}{ccc}
0 & 0 & 0 \\
4 & 4 & 0 \\
4 & 4 & 0 \\
\end{array}
\right]_{2h}^{h}
\end{equation*}
\end{minipage}
\caption{\footnotesize{ Upper, the nine points of $\mathcal{N}(x,y)$ and the green boundary of the rectangle $\mathcal{T}$. The bold red point is on the coarser and finer grids, while the little red points are on the finer grid. The arrows represent the action of the restriction operators. Below, the respective stencil to be used. }}
\label{fig:Neigh}
\end{figure}

Now consider the maximum full rectangle $\mathcal{T}$ with vertices belonging to $\mathcal{N}(x,y)$ and such that $\mathcal{T} \cap D_h \subseteq \mathcal{N}(x,y) \cap I_h$ (see Fig. \ref{fig:Neigh}). Therefore, the stencil we use in $(x,y)$ to transfer $\wh$ to a coarse grid depends on the size of $\mathcal{T}$. In fact, let $\mathcal{T} \cap D_h$ be a $3 \times 3$ points (i.e. $\mathcal{N}(x,y) \subseteq I_h$). In this case we can use the standard full-weighting stencil (\ref{Fstencil}).

Now let $\mathcal{T} \cap D_h$ be a $3 \times 2$ points. Without loss of generality, we can suppose the vertices of $\mathcal{T}$ are $(x+jh,y+ih)$, with $j \in \left\{-1,0\right\}$, $i \in \left\{ -1,1 \right\}$. In this case, the stencil we will use is:
\begin{equation}\label{Fstencil23}
\left( \tilde{I}_{2h}^h \wh \right) (x,y) = \frac{1}{16}
\left[
\begin{array}{ccc}
2 & 2 & 0 \\
4 & 4 & 0 \\
2 & 2 & 0 \\
\end{array}
\right]_{2h}^{h} (x,y),
\end{equation}
while, if $\mathcal{T}$ is a $2 \times 2$ points, with vertex $(x+jh,y+ih)$, $j,i \in\in \left\{-1,0\right\}$, the stencil will be:
\begin{equation}\label{Fstencil22}
\left( \tilde{I}_{2h}^h \wh \right) (x,y) = \frac{1}{16}
\left[
\begin{array}{ccc}
0 & 0 & 0 \\
4 & 4 & 0 \\
4 & 4 & 0 \\
\end{array}
\right]_{2h}^{h} (x,y),
\end{equation}

This three case are summarized in Fig. \ref{fig:Neigh} (where $I_h=\Omega_h$).

Finally, for all points $(x,y) \in D_{2h}-I_h$, we set $\whh(x,y) = \tilde{w}_{2h} (x,y)$.

In such a way, we can easily define the restriction operator $I_{2h}^h \colon S(D_h) \rightarrow S(D_{2h})$ as follows:
\[
I_{2h}^h \wh = \tilde{I}_{2h}^h \left( \wh, D_h-\Omega_h, \tilde{I}_{2h}^h \left( \wh, \Omega_h, 0 \right)  \right).
\]
Note that the stencils (\ref{Fstencil}), (\ref{Fstencil22}), (\ref{Fstencil23}) can be derived as tensor products of the one-dimensional restriction (\ref{restI}), (\ref{restL}), (\ref{restR}).
% in accord with the two-dimensional version of the condition (\ref{intCond}).

\subsection{Extension of the defect}\label{extension}
\begin{figure}[!hbt]%
\begin{center}
  \includegraphics[width=8cm]{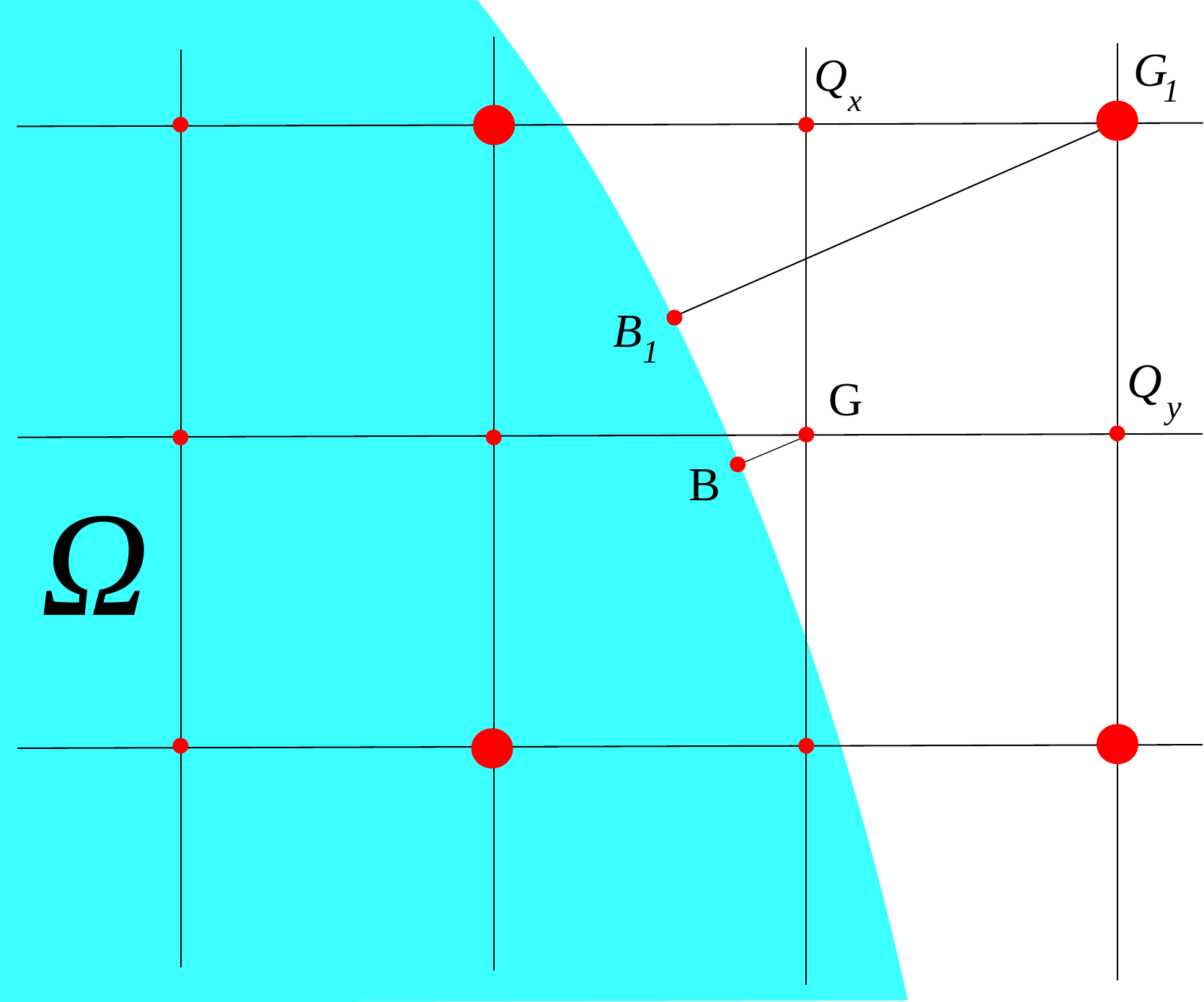}%
  \caption{\footnotesize{ Red bold and small points are grid points of $\Omega_{h}$, while red bold points are grid points of $\Omega_{2h}$. $G_1$ is a ghost point on the coarser grid, but not on the finer grid, then no value of the defect is stored in it. $Q_x$ and $Q_y$ are the two upwind near points to $G_1$. }}%
  \label{fig:extension}%
  \end{center}
\end{figure}

In every ghost point we store the defect of the boundary condition concerning that ghost point. In formulas, we have seen in step 3 of the TGCS algorithm that $\rh(G) = \left(\gh - L_h \uh \right)(G)$, for any ghost point $G$. But $\gh(G) = g(B)$ and $L_h \uh(G)$ is the reconstructed boundary condition in $B$ of the boundary operator $L$ (see (\ref{discrBC})), where $B$ is the closest boundary point to $G$ (i.e.\ the orthogonal projection on the boundary). In summary, the defect is stored in a ghost point $G$, but it is geometrically referred to a boundary point $B$ placed along the normal direction. When we switch to a coarse grid, some ghost point $G_1$ may not be ghost point in the fine grid, i.e. $\Gamma_{2h} \subseteq \Gamma_h$ is not true (see Fig.\ \ref{fig:extension}).

Then, no acceptable value of the defect is stored in $G_1$. Indeed, we expect that $\rhh$ has in the ghost point $G_1$ the defect of the boundary conditions referred to $B_1$. Hence, if we extend the defect $\rh$ outside $\Omega_h$ constant along the normal lines to the boundary, we will find $\rh\left( G_1 \right)$ as an approximation of the defect of the boundary conditions in $B_1$. After coarsening (performed using only points outside $\Omega_h$, as described before), the ghost points of the coarser grid will contain the expected values of the defect.

The extension of the defect $r_h$ is performed by solving the transport equation
\[
\frac{\partial r}{\partial \tau} + \nabla r \cdot \vec{n} = 0
\]
in a few steps of a fictitious time $\tau$, where $\vec{n} \equiv (n_x, n_y) = \nabla \phi / |\nabla \phi|$ is the unit normal vector to the level-set, while $r(x,y)$ is a continuous version of $\rh$ (i.e. $r(x,y)$ is a continuous function defined in $D-\Omega$ and such that $r(G) = \rh(G)$ for any ghost point $G$). In details, we compute few steps of the following iteration scheme:
\begin{equation}\label{exp}
\rh^{(m+1)}(P) = \rh^{(m)}(P) + \frac{\Delta \tau}{h} \left( \left( \rh^{(m)}(P) - \rh^{(m)}(Q_x) \right) \left| n_x \right| +  \left( \rh^{(m)}(P) - \rh^{(m)}(Q_y) \right) \left| n_y \right| \right)
\end{equation}
for all $P \in D_h - \left( \Omega_h \cup \Gamma_h \right)$, where $Q_x$ and $Q_y$ are the two upwind near points to $P$, i.e.
\[
Q_x=P-\sgn(n_x) \: h \: \vec{i}, \; \; \; Q_y=P-\sgn(n_y) \: h \: \vec{j}.
\]
However, it is sufficient to perform the iteration (\ref{exp}) only in a narrow band with width $3 \, h$.
In order to speed up the extension, we can perform (\ref{exp}) in a Gauss-Seidel fashion, sorting points in $D_h - \left( \Omega_h \cup \Gamma_h \right)$ by the distance from the boundary (it can be done using the distance function $\phi$), and starting the computation $\rh^{(m+1)}(P)$ in (\ref{exp}) from the closest ghost point $P$ to $\Gamma$.
%up to the ghost point in the coarse grid with the largest distance from $\Gamma$.

\subsubsection{Interpolation}
%se hai modificato i'interpolation operator in 1D, modificalo anche qua
Since the interpolation operator acts on the error, which is continuous across the boundary, we just use the standard bilinear interpolation operator:
\begin{equation*}
I_{h}^{2h}=
\frac{1}{4}
\left]
\begin{array}{ccc}
1 & 2 & 1 \\
2 & 4 & 2 \\
1 & 2 & 1 \\
\end{array}
\right[_{h}^{2h}.
\end{equation*}

\section{Numerical tests}\label{numtests}
In all the following numerical tests we always choose the Dirichlet and Neumann parts of $\Gamma=\partial \Omega$ as:
\[
\Gamma_D=\left\{ (x,y) \in \Gamma \colon x \leq 0 \right\}, \; \; \; \Gamma_N=\left\{ (x,y) \in \Gamma \colon x > 0 \right\}.
\]
The Local Fourier Analysis (LFA) is a powerful tool to obtain the theoretically convergence factor by analyzing separately the action of different parts of the multigrid algorithm to high and low frequency components of the error. For a detailed explanation of the LFA, we refer to~\cite[Ch. 4]{Trottemberg:MG}. 

Before to apply the LFA, one has to be sure that the relaxation operator (\ref{relaxO}) has the \textit{smoothing property}. Roughly speaking, the \textit{smoothing property} is the property to dump high frequency components of the error, in order to make it smooth after few relaxation sweeps.

When the multigrid algorithm applies to a regular rectangular domain, the LFA and smoothing analysis are well studied. In the case of arbitrary domain, as Achi Brandt points out in~\cite[pag. 587]{Trottemberg:MG}, there are some boundary related difficulties about the discretization and relaxation near the boundary: 
\begin{itemize}
\item There is no a general smoothing analysis when the boundary is not aligned with the grid;
\item The residuals should be reduced near boundaries more than in the interior;
\item The coarsest grid has not to be too coarse, because it should catch the curvature of the boundary in order to guarantee the convergence.
\end{itemize}
Now, we perform a numerical test in order to check if the convergence factor is close to the predicted one by LFA, which is obtained for rectangular domain with periodic boundary conditions, i.e.\ without taking into account boundary effects.
Note that the multigrid algorithm described before may be seen as an iterative scheme:
\[
\uh^{(m+1)} = M_h \: \uh^{(m)} + \vec{b}_h
\]
for some matrix $M_h$ and vector $\vec{b}_h$. We call $\rho$ the convergence factor, which is the spectral radius of the matrix $M_h$. For rectangular domain with periodic boundary conditions and constant coefficients, the convergence factor is said to be local and it is denoted by $\rho_{loc}$.
The convergence factors predicted by LFA for Gauss-Seidel LEX relaxation and FW restriction operator are listed in Table \ref{table:rholoc} (see~\cite[pag. 117]{Trottemberg:MG}).

\begin{table}[!hbt]
\captionsetup{width=0.79\textwidth}
\caption{ \footnotesize{ Predicted convergence factor $\rho_{loc}$ by LFA for GS-LEX and FW restriction operator.}} % title of Table
\centering      % used for centering table
\begin{tabular}{||c || c | c | c | c ||}  % centered columns
\hline\hline                        %inserts double horizontal lines
$\nu=\nu_1+\nu_2$ & 1 & 2 & 3 & 4 \\ [0.5ex] % inserts table
\hline                    % inserts single horizontal line
$\rho_{loc}$ & 0.400 & 0.193 & 0.119 & 0.084 \\
%heading
\hline     %inserts single line 
 \end{tabular} 
 \label{table:rholoc}
 \end{table} 
 
In all the numerical tests we perform, the convergence factor is estimated as the ratio of consecutive defects, i.e.:
\[
\rho=\rho^{(m)} = \frac{\left\| \rh^{(m)} \right\|_{\infty}}{\left\| \rh^{(m-1)} \right\|_{\infty}}
\]
for $m$ very large.
In order to avoid difficulties related to numerical instability related to the machine precision, we will always use the homogeneous model problem as a test, namely (\ref{mainpr}) with $f=g_D=g_N=0$, and perform the multigrid algorithm starting from an initial guess different from zero. Since we are just interested at the convergence factor and not at the numerical solution itself (which approaches zero indefinitely for homogeneous problem), a reasonable stopping criterion will be
\[
\frac{ \left|\rho^{(m)}-\rho^{(m-1)} \right| }{ \rho^{(m)} } < 10^{-3}.
\]
Note that, since we want to study the efficiency of the multigrid components proposed in this paper (smoother, restriction, ...), it is sufficient to study basic kind of multigrid such as V-cycle and W-cycle, while a more efficient algorithm (such as FMG) can be easily derived.

\subsection{1D numerical test}\label{ex:1D}
Referring to Sec. \ref{1Dsection}, let us choose $[a,b]=[-0.743,0.843] \subseteq [-1,1]$. The finest grid is obtained dividing the whole computational domain $[-1,1]$ into $N=64$ subintervals, while the coarsest grid is obtained dividing $[-1,1]$ into $Nc=8$ subintervals. The computed convergence factors are very close to those ones predicted by LFA (Table \ref{table:rholoc}), namely $\rho=0.185$ for $\nu=2$ and $\rho=0.122$ for $\nu=3$.

\subsection{An initial test in 2D}
\begin{figure}[!hbt]%
\captionsetup{width=0.79\textwidth}
  \begin{center}
  \includegraphics[width=0.60\textwidth]{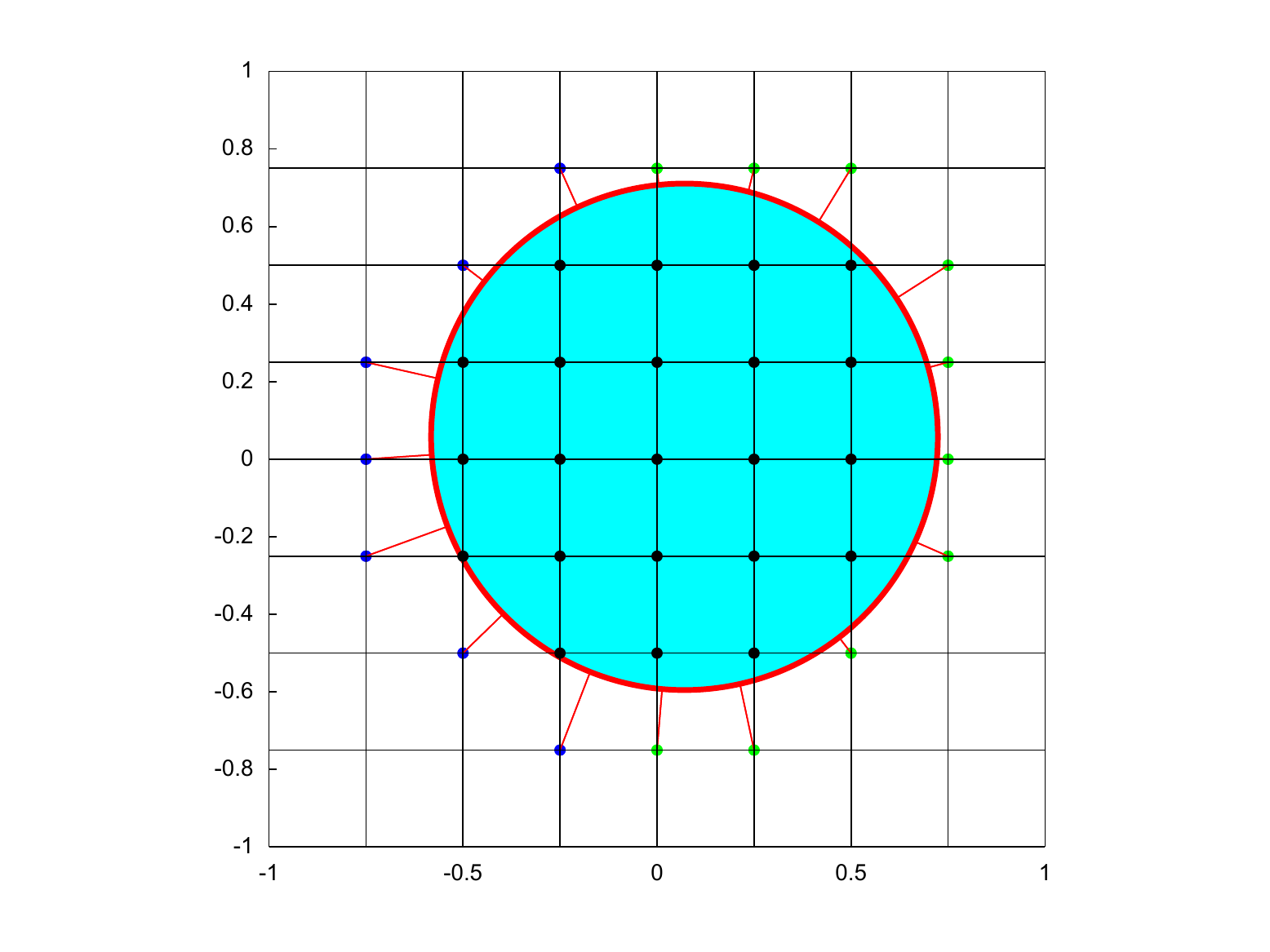}%
  \caption{\footnotesize{ Circular domain and the coarsest grid used to capture the curvature. Blue ghost points refer to Dirichlet condition, while green ghost points refer to Neumann condition. Red lines are normal to the boundary.}}
   \label{fig:circle1}
  \end{center}
\end{figure}

We start testing the multigrid algorithm on a circular domain $\Omega$ with center $(\sqrt{2}/20,\sqrt{3}/30)$ and radius $r=0.563$ (Fig. \ref{fig:circle1}).

The measured convergence factors for TGCS, $V$-cycle and $W$-cycle are listed in Table \ref{table:badrho}.

 \begin{table}[!hbt]
\captionsetup{width=0.79\textwidth}
\caption{ \footnotesize{ Measured convergence factor $\rho$ with $\nu_1=\nu_2=1$ on the left and with $\nu_1=2$, $\nu_2=1$ on the right.}} % title of Table
\begin{minipage}{0.45\textwidth}
\centering      % used for centering table
\begin{tabular}{||c || c | c | c ||}  % centered columns (4 columns)
\hline\hline                        %inserts double horizontal lines
$N$ & TGCS & $V$-cycle & $W$-cycle \\ [0.5ex] % inserts table
%heading
\hline                    % inserts single horizontal line
64 & 0.67 & 0.68 & 0.71 \\
128 & 0.68 & 0.73 & 0.68 \\
256 & 0.70 & 0.71 & 0.70 \\
\hline     %inserts single line 
 \end{tabular} 
\end{minipage} 
\begin{minipage}{0.45\textwidth}
\centering      % used for centering table
\begin{tabular}{||c || c | c | c ||}  % centered columns (4 columns)
\hline\hline                        %inserts double horizontal lines
$N$ & TGCS & $V$-cycle & $W$-cycle \\ [0.5ex] % inserts table
%heading
\hline                    % inserts single horizontal line
64 & 0.58 & 0.72 & 0.58 \\
128 & 0.58 & 0.73 & 0.59 \\
256 & 0.61 & 0.83 & 0.60 \\
\hline     %inserts single line 
 \end{tabular} 
 \end{minipage}
  \label{table:badrho}
 \end{table} 
 
As we can see, the measured convergence factor is far from the predicted one by LFA (Table \ref{table:rholoc}). Then, some boundary effect degrades the convergence factor. Note that in 1D such boundary effects do not degrade the convergence factor (Ex. \ref{ex:1D}), because we have only two boundary points, and the degradation is due to the oscillating behavior of the residual on the tangential direction to the boundary, that does not exist in 1D. Then, in 2D we must smooth the error also along the tangential direction to the boundary.
To overcome this difficulty, we apply, after a single relaxation and at each grid level, $\lambda$ extra relaxation sweeps on all ghost points $\Gamma_h$ and on all inside grid points of $\Omega_h$ within $\delta>0$ distance from the boundary (the extra computational work is $O(N)$, then negligible as $N \rightarrow \infty$).
It can be proved numerically that a good choice of these parameters will be:
\[
\lambda=5, \; \; \; \delta=3 \: h.
\]
The explanation of the optimal value $\lambda=5$ is the following: the degradation observed in Table \ref{table:badrho} is an indication that the error decays much slower at the boundary. Assuming that the convergence factor in Table \ref{table:badrho} is essentially the convergence factor at the boundary, $\rho_B$, we want to match it with the convergence factor at the bulk, therefore $\lambda_{opt}$ is the smallest value of $\lambda$ for which $\rho_B^{\lambda+1} \leq \rho_I$. The value $\rho_I$, in turn, can be computed as the convergence factor for large value of $\lambda$.

\begin{figure}[!hbt]%
\begin{minipage}{0.49\textwidth}
  \begin{center}
  \includegraphics[width=0.75\textwidth]{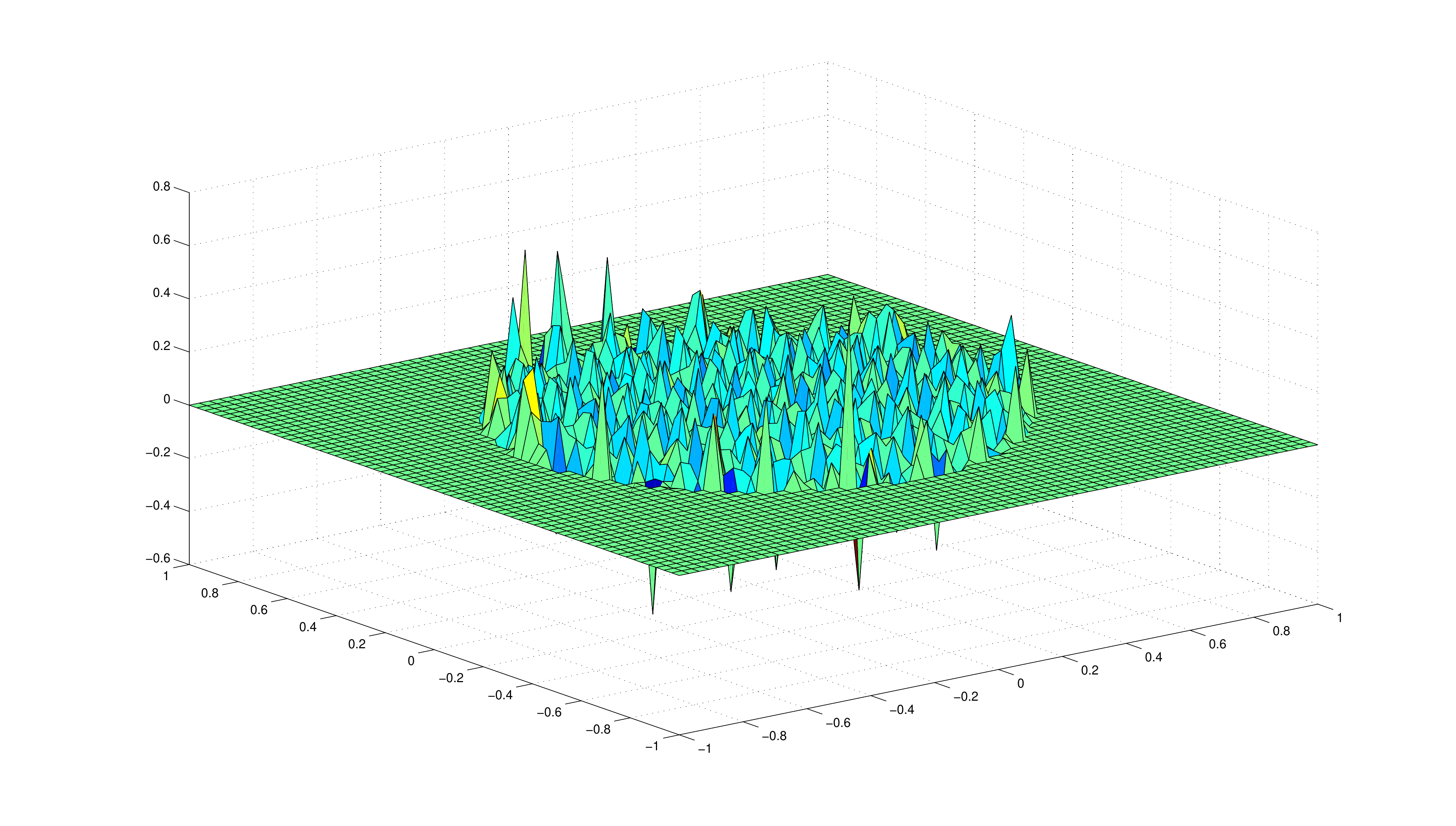}%
  \end{center}
\end{minipage}
\begin{minipage}{0.49\textwidth}
  \begin{center}
  \includegraphics[width=0.75\textwidth]{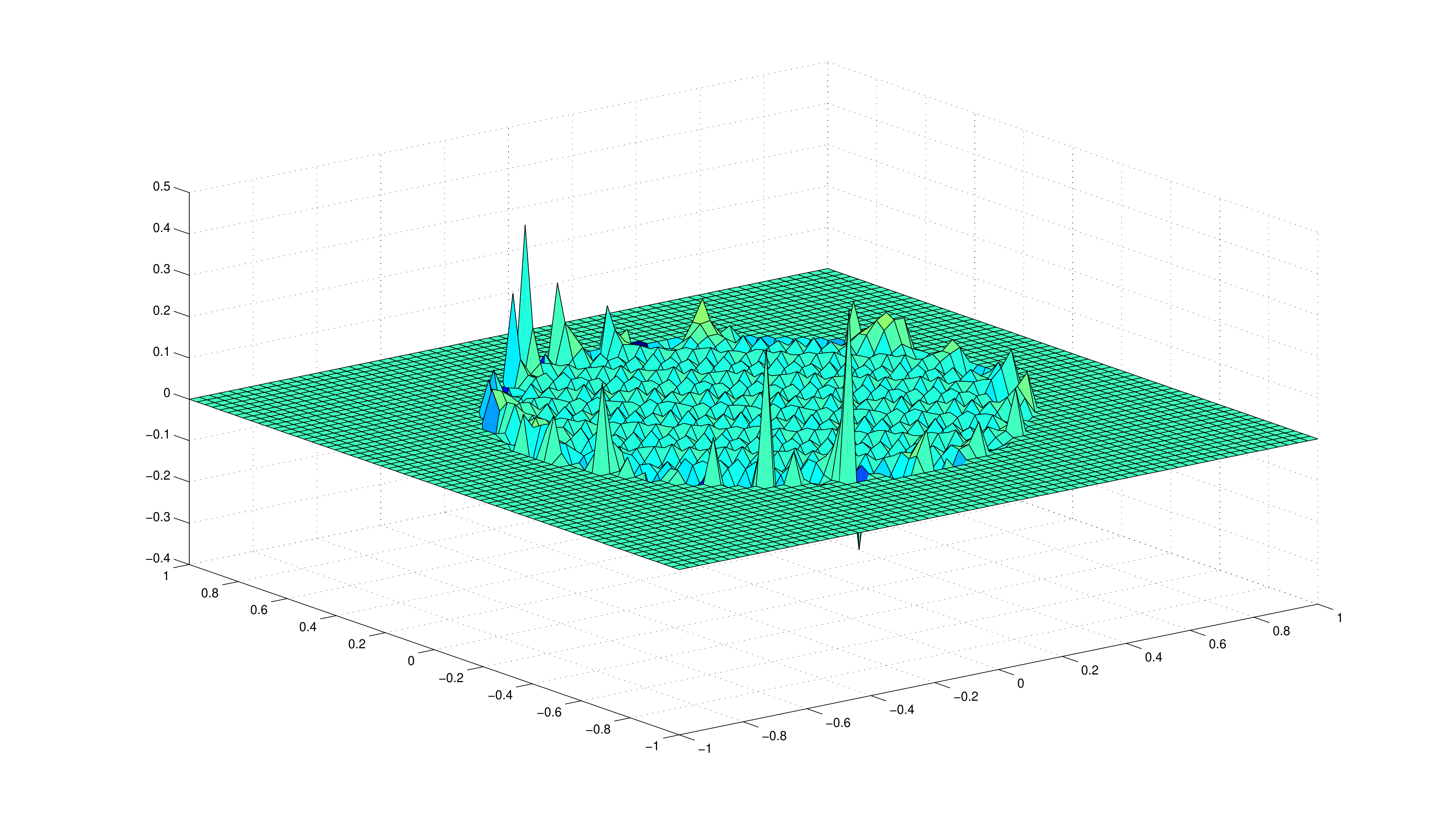}%
  \end{center}
\end{minipage}
\begin{minipage}{0.49\textwidth}
  \begin{center}
  \includegraphics[width=0.75\textwidth]{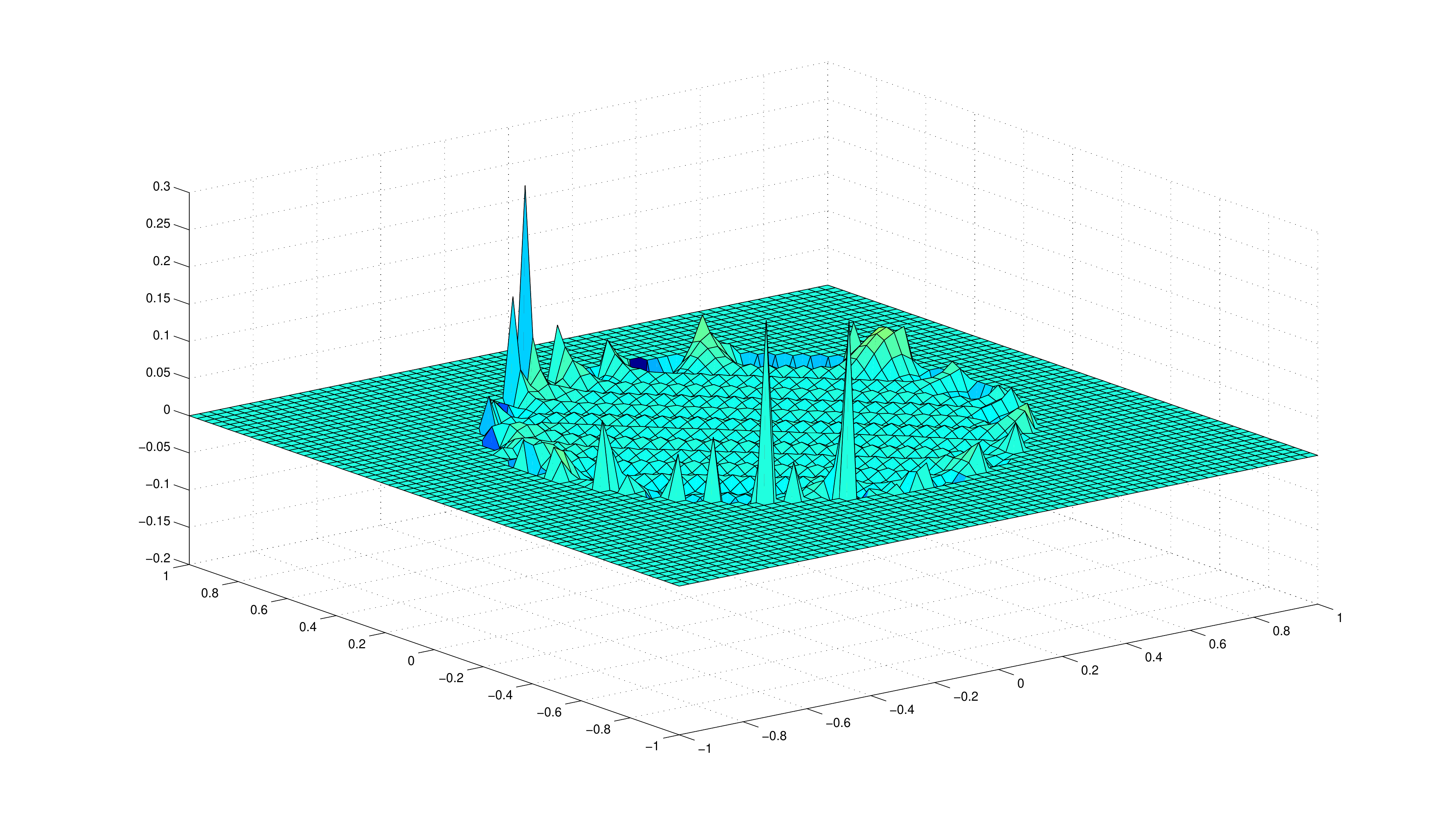}%
  \end{center}
\end{minipage}
\begin{minipage}{0.49\textwidth}
  \begin{center}
  \includegraphics[width=0.75\textwidth]{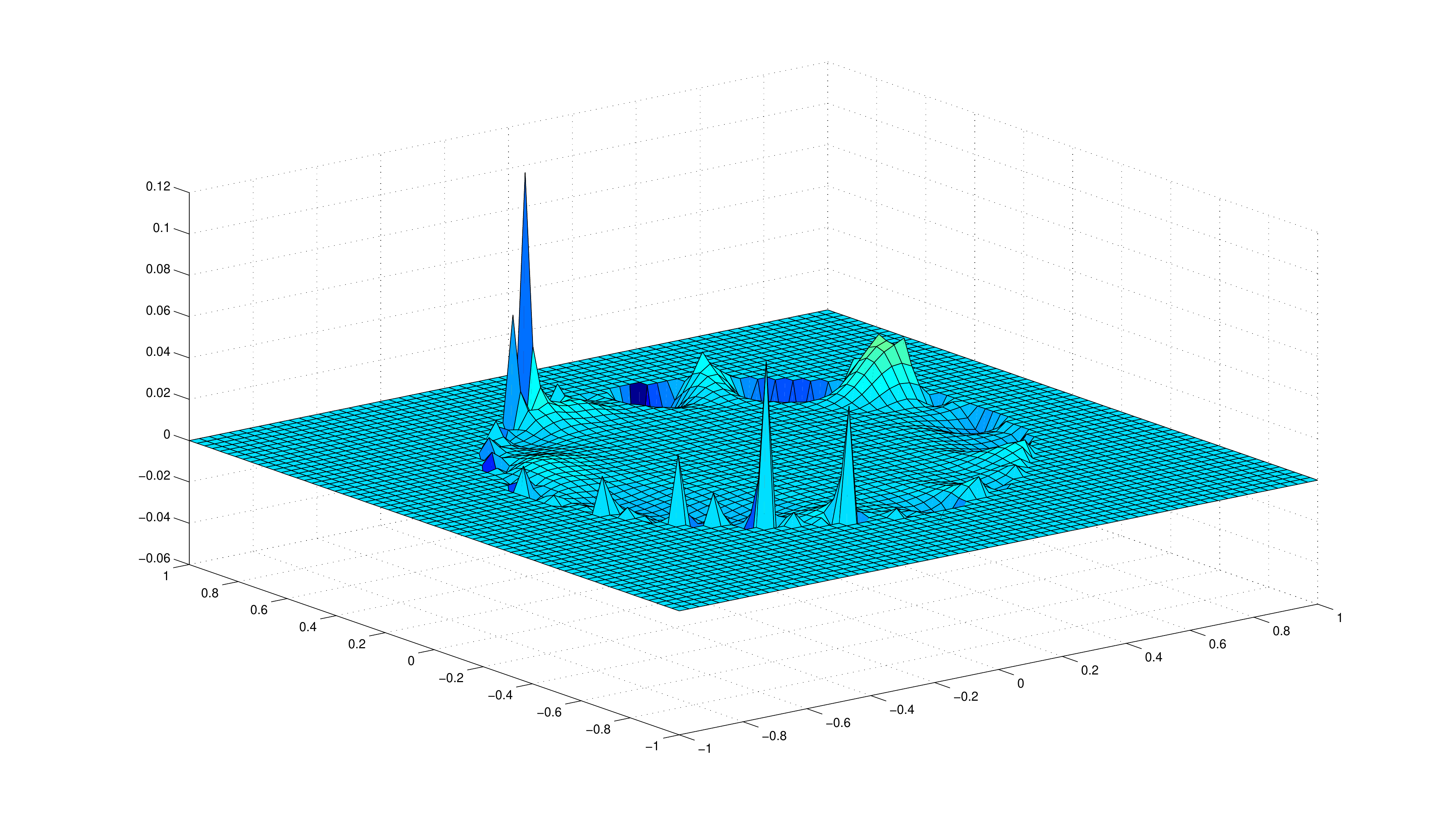}%
  \end{center}
\end{minipage}
\caption{\footnotesize{ High frequency initial error after $1$ (up-left), $3$ (up-right), $5$ (down-left), $10$ (down-right) relaxation sweeps and $\lambda=0$ extra-relaxations. }}
\label{fig:extraNO}
\vskip 0.3 cm 
\begin{minipage}{0.49\textwidth}
  \begin{center}
  \includegraphics[width=0.75\textwidth]{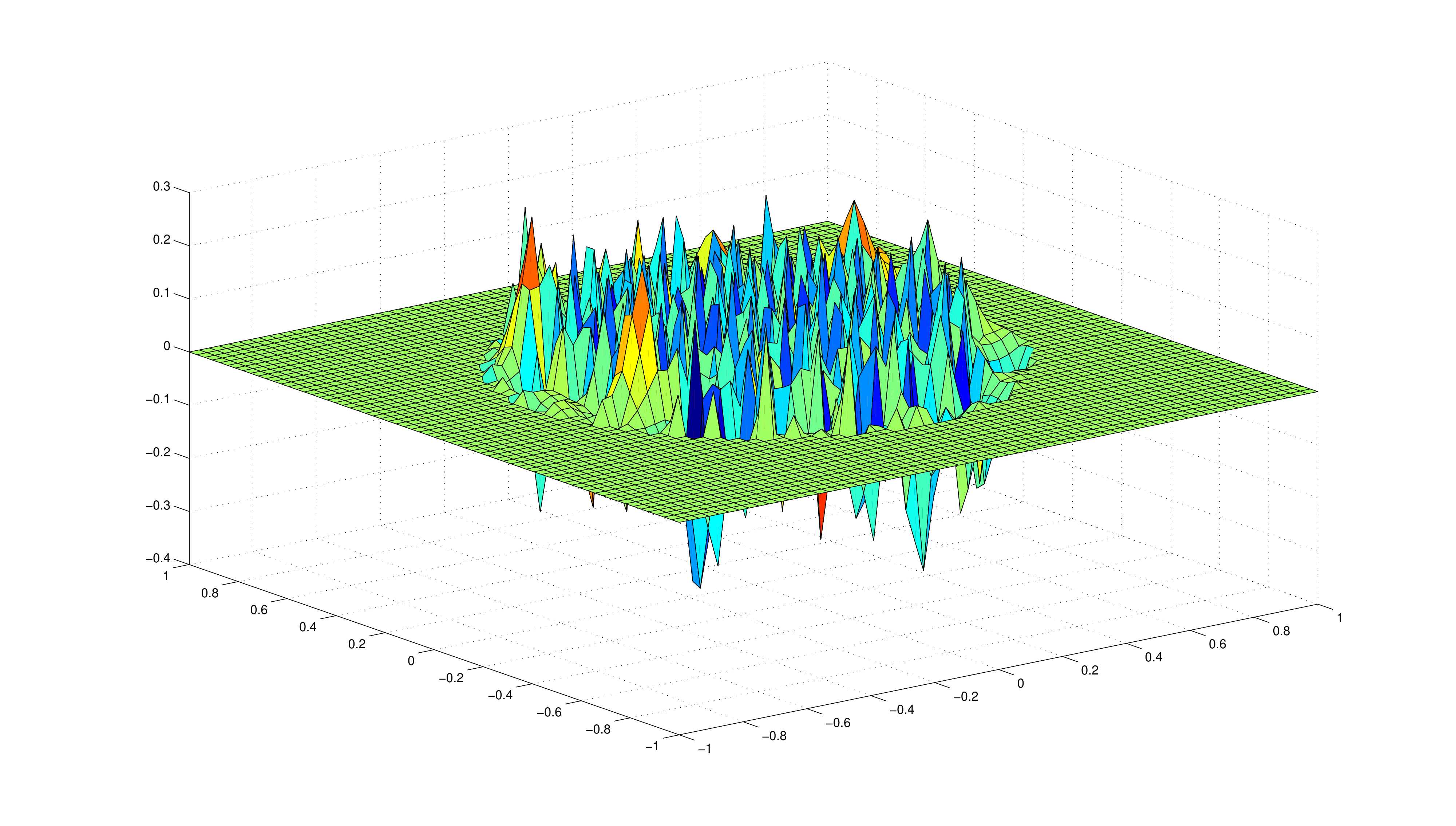}%
  \end{center}
\end{minipage}
\begin{minipage}{0.49\textwidth}
  \begin{center}
  \includegraphics[width=0.75\textwidth]{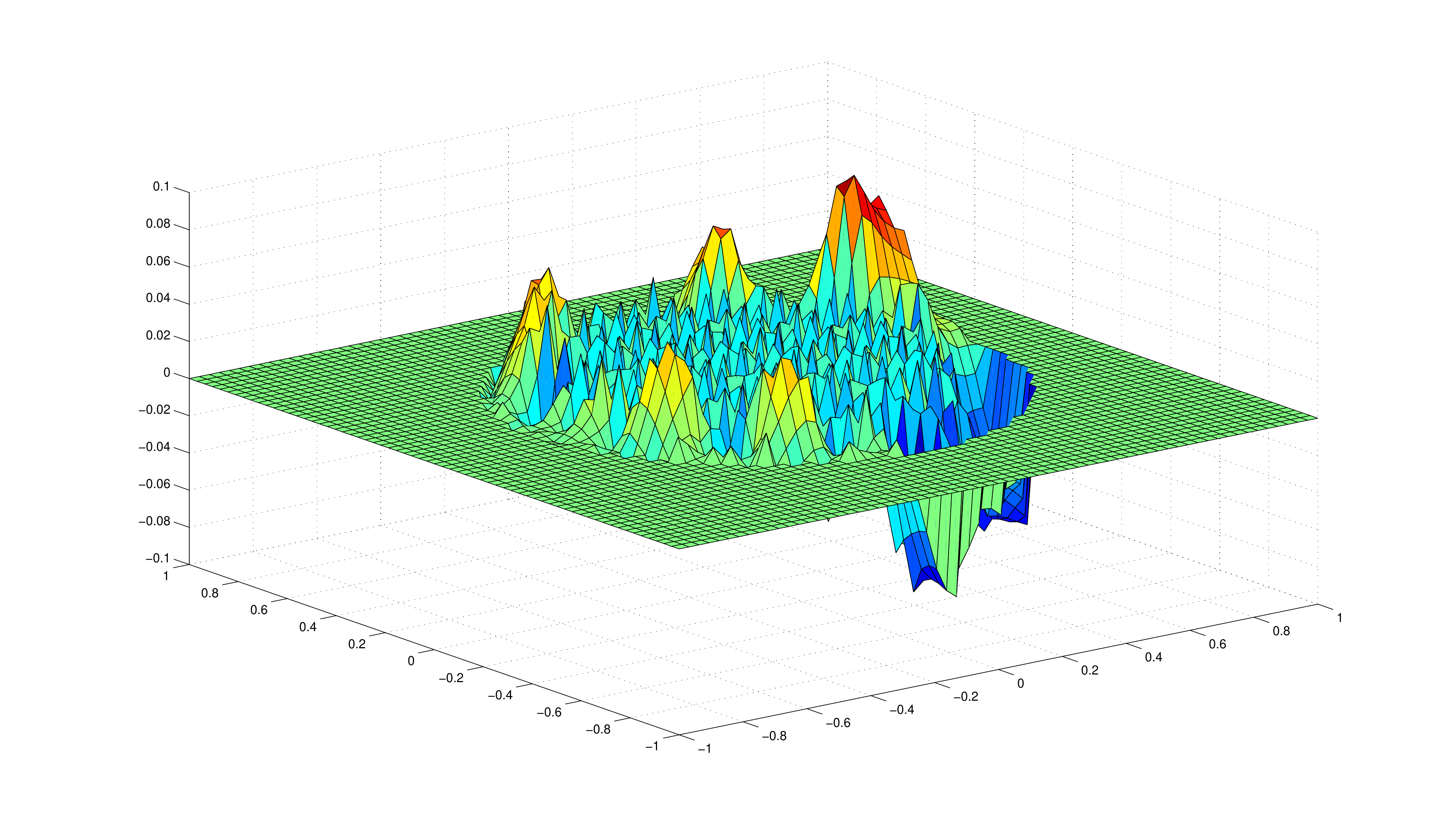}%
  \end{center}
\end{minipage}
\begin{minipage}{0.49\textwidth}
  \begin{center}
  \includegraphics[width=0.75\textwidth]{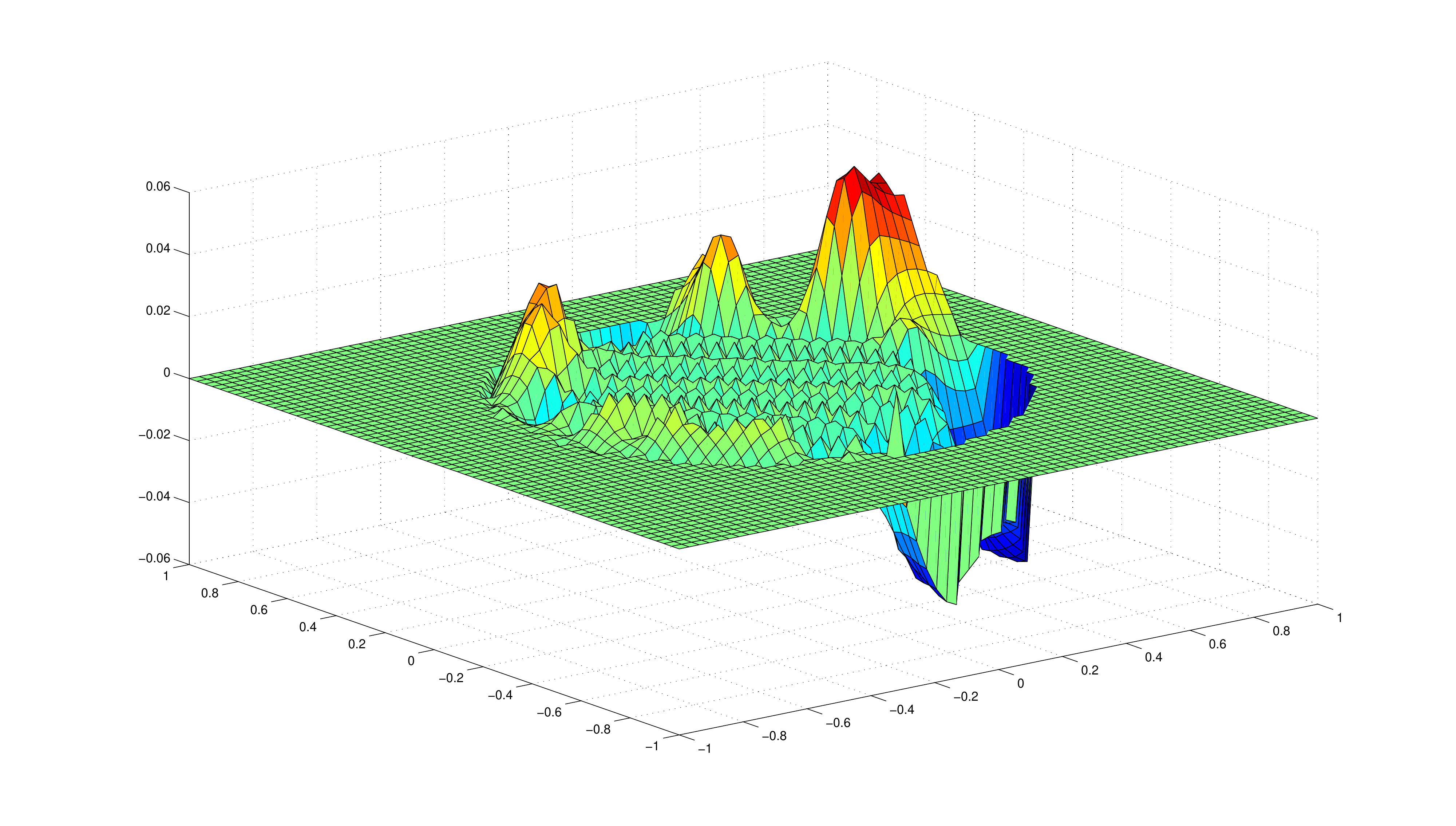}%
  \end{center}
\end{minipage}
\begin{minipage}{0.49\textwidth}
  \begin{center}
  \includegraphics[width=0.75\textwidth]{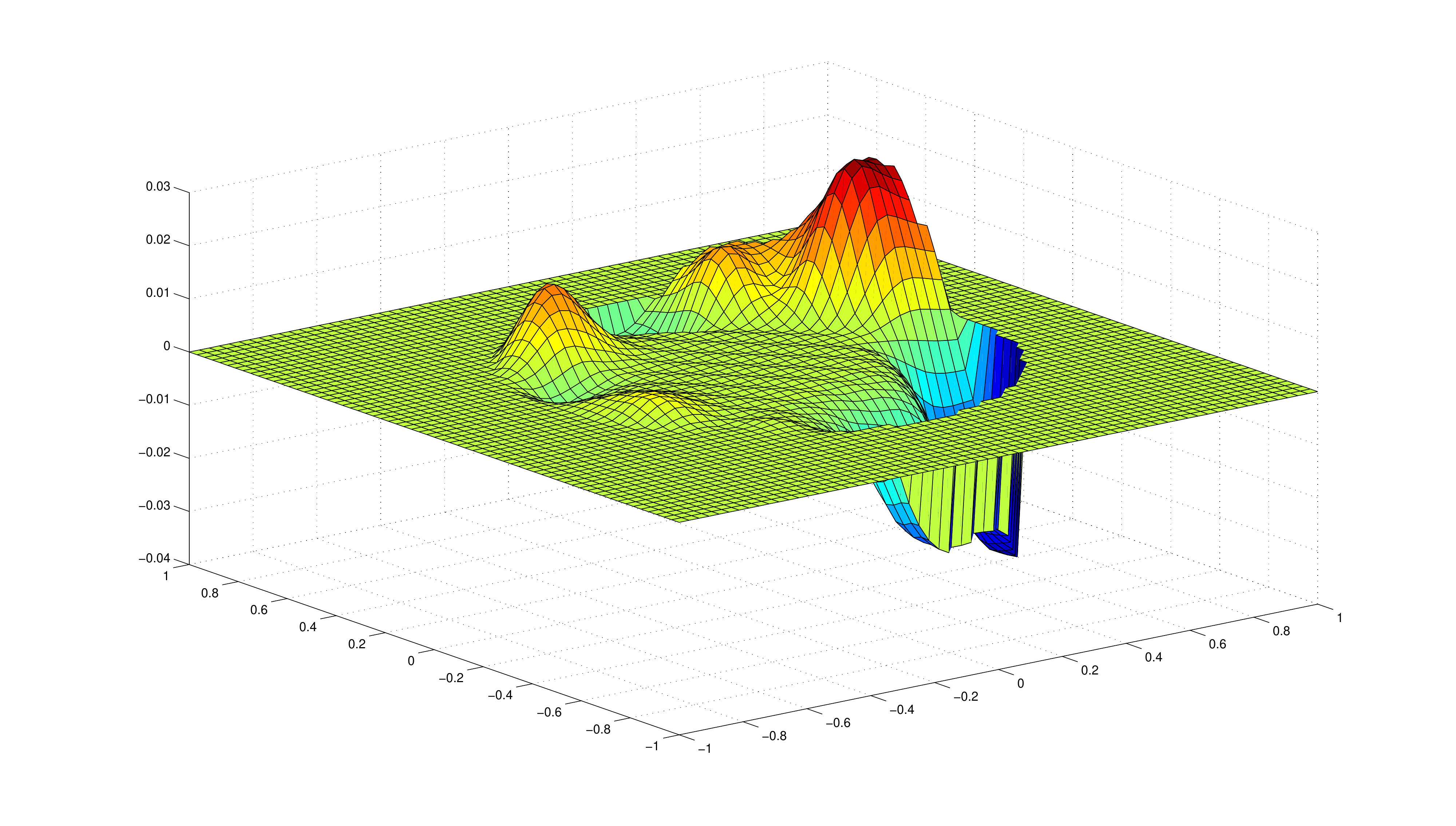}%
  \end{center}
\end{minipage}
\caption{\footnotesize{ High frequency initial error after $1$ (up-left), $3$ (up-right), $5$ (down-left), $10$ (down-right) relaxation sweeps and $\lambda=5$ extra-relaxations. }}
\label{fig:extraYES}
\end{figure}

Investigating the smoothing property, we observe that choosing the initial error as an high frequency component, the error is not smoothed after few relaxation sweeps. While, if we add the extra-relaxations, the error become sufficiently smooth (Figs. \ref{fig:extraNO}-\ref{fig:extraYES}).

\subsection{Some numerical results}
In this section we confirm numerically the improvement of the convergence factor if we apply extra-relaxations, and we compare the relaxations with other well-knowns alternative such as the Kaczmarz and the block relaxations.
In all numerical tests, we choose an arbitrary domain $\Omega$ assigning a level-set function $\phi_0$. Then we reinitialize it by the procedure described in Section \ref{levelset}, obtaining the signed distance function $\phi$. Afterwards, we perform the multigrid technique applying the $W$-cycle algorithm instead of the $V$-cycle, to ensure the independence of the convergence factor $\rho$ from the step size $h$ (as explained for example in~\cite[pag. 78]{Trottemberg:MG}). Several tests are performed for each domain, based on the different size of the finest and coarsest grids. The finest grid is obtained dividing the whole computational domain $D$ into $N$ subintervals in each Cartesian direction, while the coarsest grid is obtained replacing $N$ with $N_c$. The solution on the coarsest grid is obtained by a direct solver.

\subsubsection{Circular domain}\label{NumTests:circle}
\begin{figure}[!hbt]%
\begin{minipage}{0.49\textwidth}
  \begin{center}
  \includegraphics[width=0.95\textwidth]{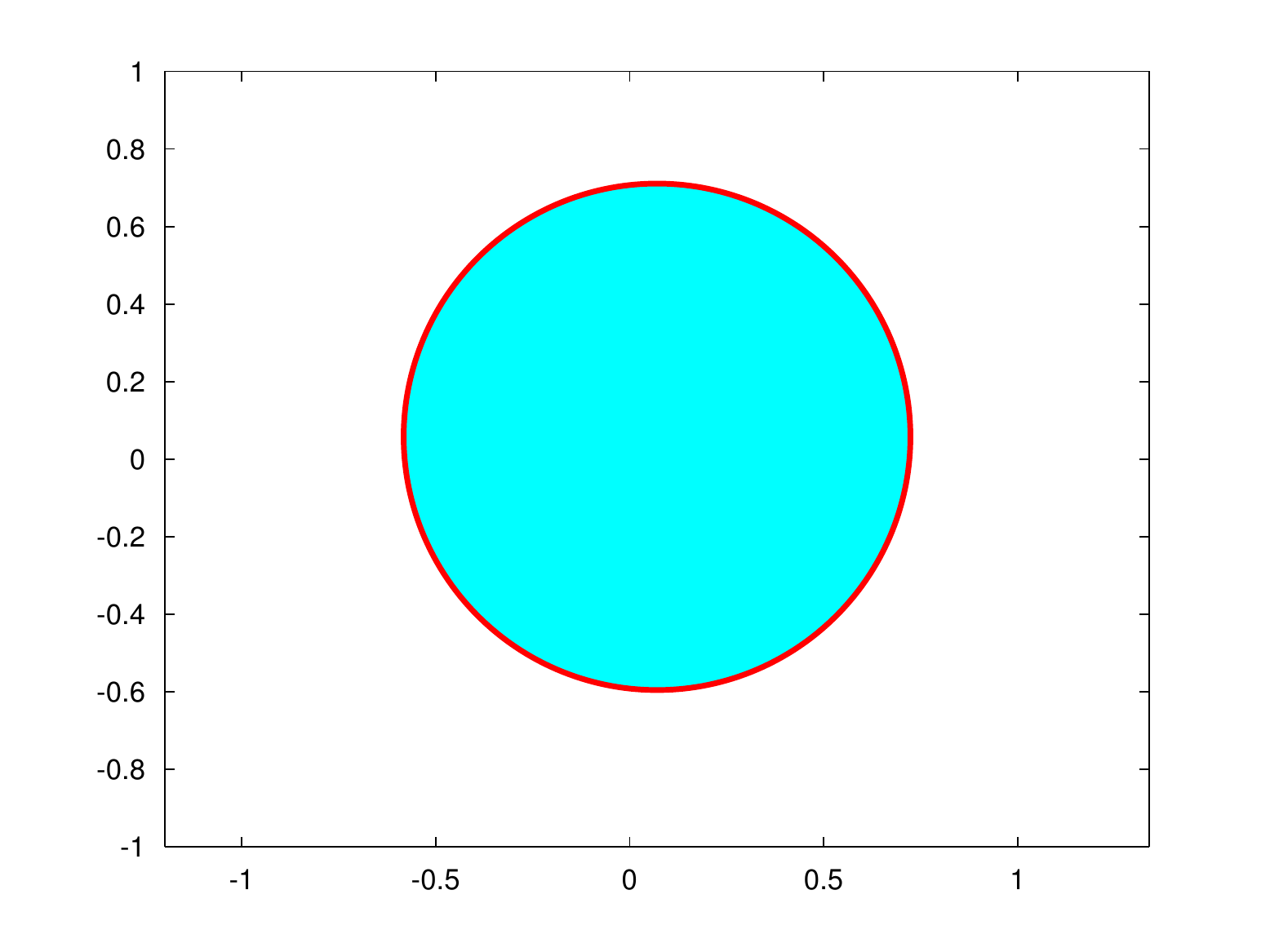}%
  \end{center}
\end{minipage}
\begin{minipage}{0.49\textwidth}
  \begin{center}
  \includegraphics[width=0.95\textwidth]{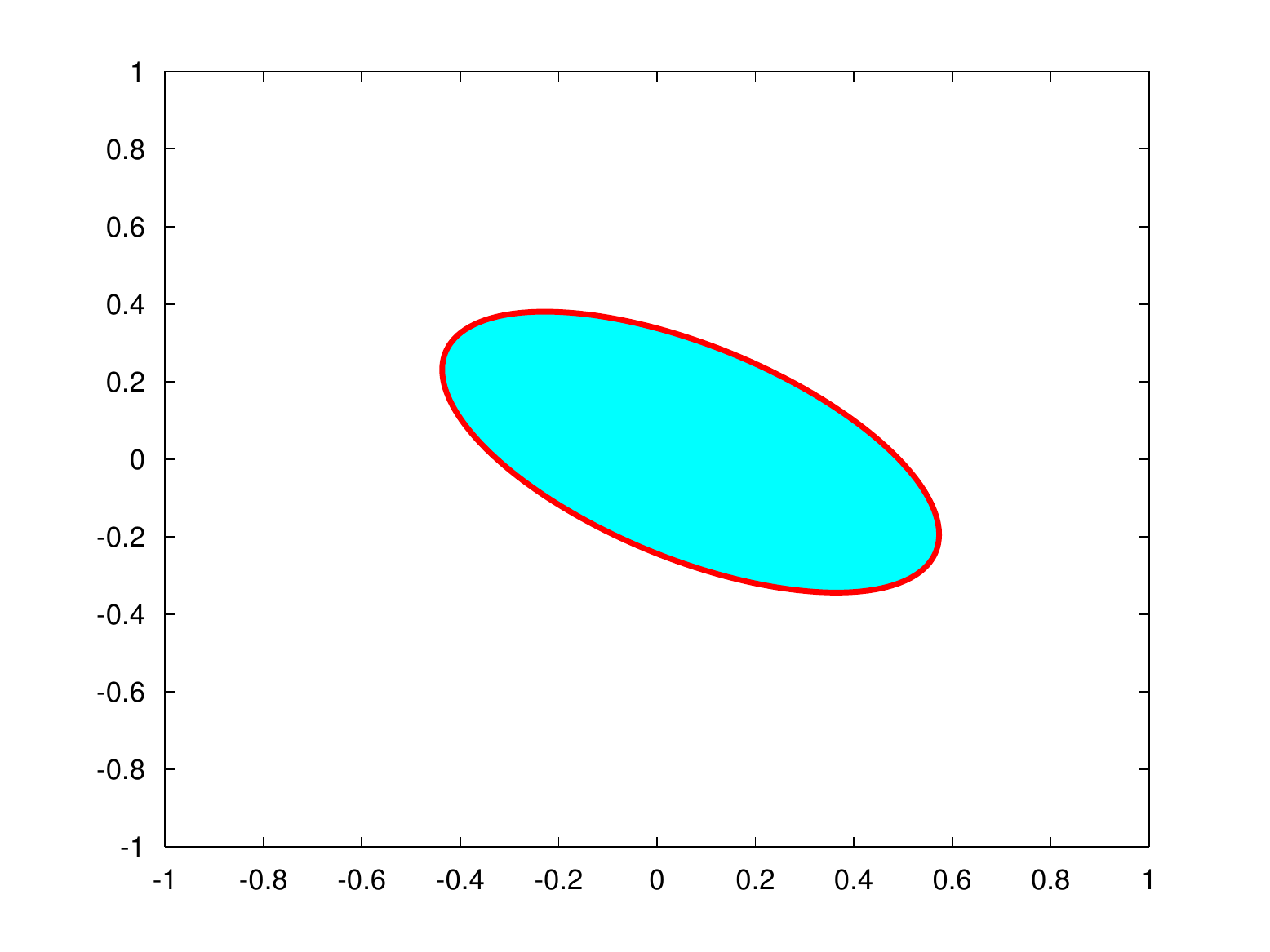}%
  \end{center}
\end{minipage}
\begin{minipage}{0.49\textwidth}
  \begin{center}
  \includegraphics[width=0.95\textwidth]{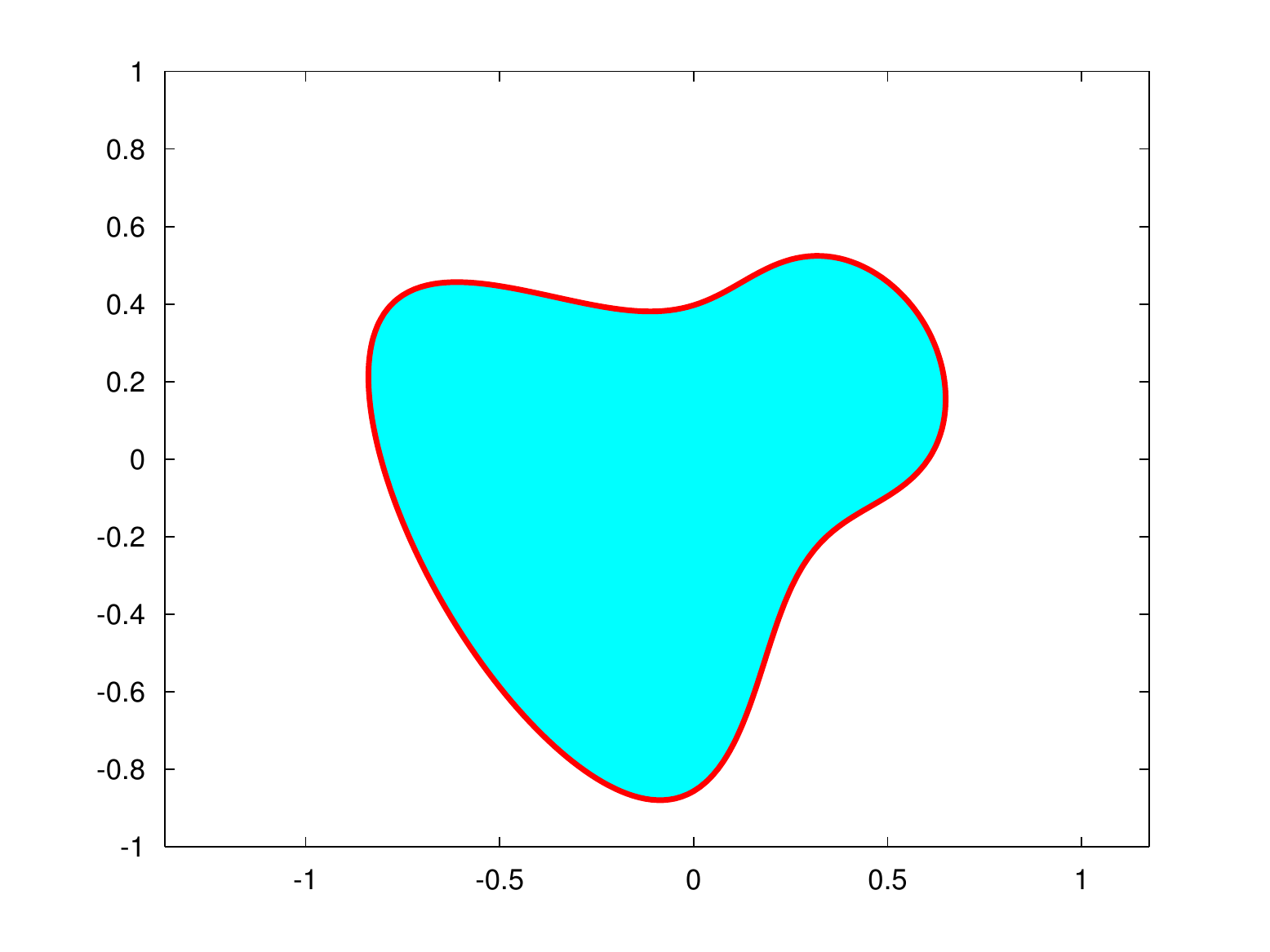}%
  \end{center}
\end{minipage}
\begin{minipage}{0.49\textwidth}
  \begin{center}
  \includegraphics[width=0.95\textwidth]{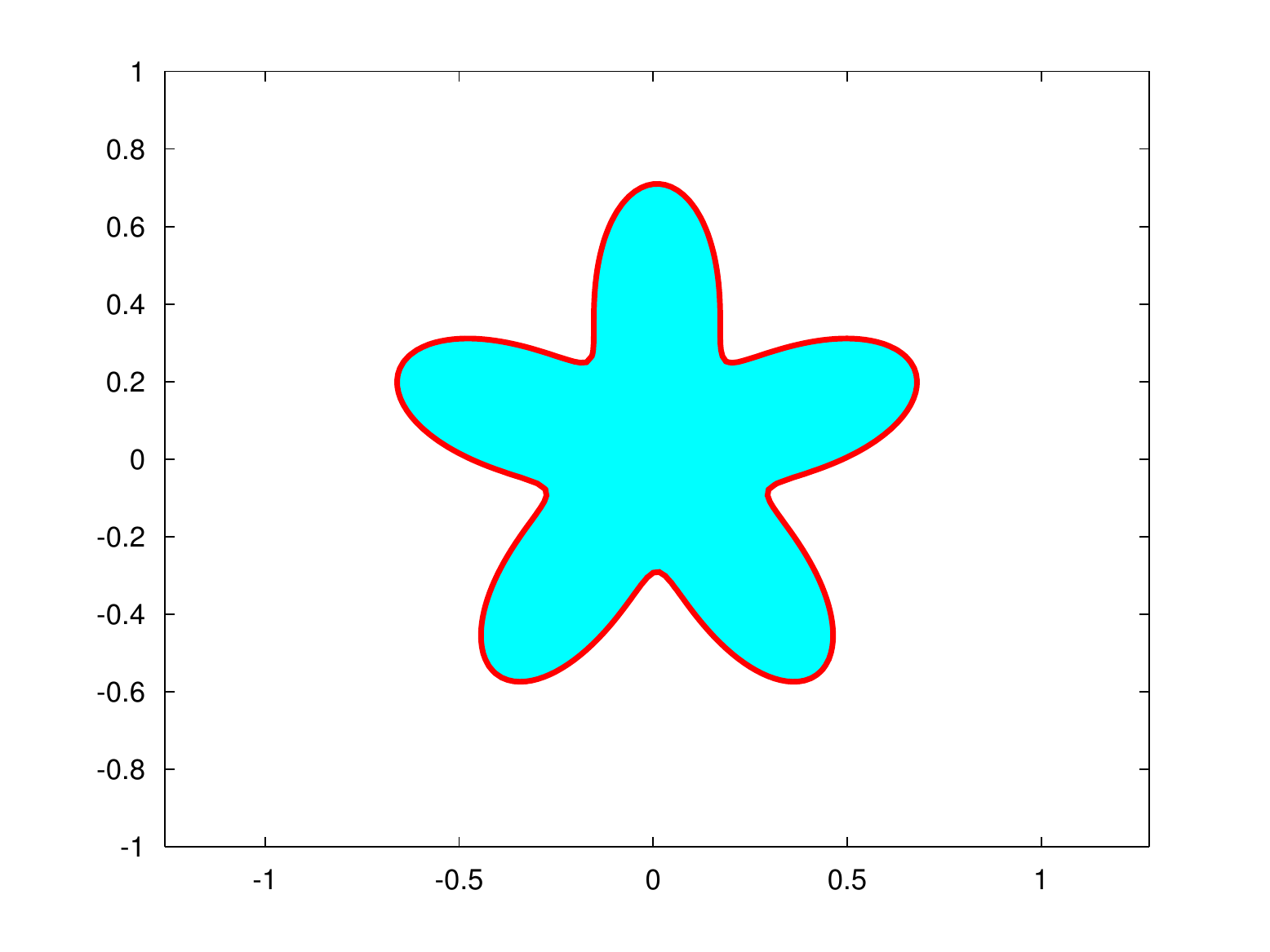}%
  \end{center}
\end{minipage}
\caption{\footnotesize{ Different domains used in the numerical tests: Example \ref{NumTests:circle} (top-left), \ref{NumTests:ellipse} (top-right), \ref{NumTests:saddle} (bottom-left), \ref{NumTests:flower} (bottom-right).}}
\label{fig:domains}
\end{figure}

In this case we can choose as a level-set function directly the signed distance function, which is known analytically:
\[
\phi(x,y) = \sqrt{(x-\sqrt{2}/20)^2 + (y-\sqrt{3}/30)^2} - 0.563.
\]
The zero level-set is represented in Fig \ref{fig:domains} (top-left).
Different value of the convergence factor are listed in Table \ref{table:rhoC} (for $\nu=\nu_1+\nu_2=2$ and $\nu=3$). They are really improved with respect to those obtained without extra-relaxations (Table \ref{table:badrho}) for the same test.

 \begin{table}[!hbt]
\captionsetup{width=0.79\textwidth}
\caption{ \footnotesize{Convergence factor of the numerical test of Section \ref{NumTests:circle}. We use $N \times N$ number of grid points in the finest grid; $N_c \times N_c$ number of grid points in the coarsest grid. Left table: $\nu=\nu_1+\nu_2=2$, right table $\nu=\nu_1+\nu_2=3$. }} % title of Table
\begin{minipage}{0.45\textwidth}
\centering      % used for centering table
\begin{tabular}{|| c c || c | c | c | c | c ||}  % centered columns (4 columns)
\hline\hline                        %inserts double horizontal lines
 & $N$ & 16 & 32 & 64 & 128 & 256  \\ [0.5ex] % inserts table
$N_c$ & &  &  &  &  & \\ 
\hline\hline                        %inserts double horizontal lines
8 & & 0.052 & 0.053 & 0.11 & 0.13 & 0.14\\ 
16 & &  & 0.061 & 0.11 & 0.13 & 0.14\\ 
32 & &  &  & 0.11 & 0.13 & 0.14\\ 
64 & &  &  &  & 0.13 & 0.14\\ 
128 & &  &  &  &  & 0.14\\ 
\hline\hline                        %inserts double horizontal lines
 \end{tabular} 
\end{minipage} 
\begin{minipage}{0.45\textwidth}
\centering      % used for centering table
\begin{tabular}{|| c c || c | c | c | c | c ||}  % centered columns (4 columns)
\hline\hline                        %inserts double horizontal lines
 & $N$ & 16 & 32 & 64 & 128 & 256  \\ [0.5ex] % inserts table
$N_c$ & &  &  &  &  & \\ 
\hline\hline                        %inserts double horizontal lines
8 & & 0.06 & 0.03 & 0.09 & 0.08 & 0.08\\ 
16 & &  & 0.04 & 0.09 & 0.08 & 0.08\\ 
32 & &  &  & 0.09 & 0.08 & 0.08\\ 
64 & &  &  &  & 0.09 & 0.08\\ 
128 & &  &  &  &  & 0.09\\ 
\hline\hline                        %inserts double horizontal lines
 \end{tabular} 
 \end{minipage}
  \label{table:rhoC}
 \end{table} 

\subsubsection{Comparison with the Kaczmarz and the block relaxations}
Note that the relaxation scheme (\ref{jac}), (\ref{left2}), (\ref{right2}) is composed by a Gauss-Seidel iteration over inner grid points and a suitable relaxation over ghost points (boundary conditions). As an alternative to the relaxation of the boundary condition, we can use the Kaczmarz relaxation~\cite{Kaczmarz:Kiteration} near the boundary, which is known to be unconditionally convergent. Let us recall the Kaczmarz iteration scheme for a subset of equations $\mathcal{J} \subseteq \left\{\, \ldots,N_i+N_g \right\}$ of a linear system $Lu=f$:
\begin{equation*}
\boxed{
\begin{split}
 & u^{TEMP}=u^{(m)}, \\
 & \mbox{for } j \in \mathcal{J} \mbox{ do: } \; \; \;
  u^{TEMP} \colon \! \! = u^{TEMP} + \frac{f_j-<l_j,u^{TEMP}>}{\left\| l_j \right\|_2^2} l_j^T, \\
 & u^{(m+1)} = u^{TEMP}.
\end{split}
}
\end{equation*}
The symbol $<\cdot,\cdot>$ denotes the inner product operator and $l_j$ is the $j$-th row of the matrix $L$.
If we choose $J=\left\{\, \ldots,N_i+N_g \right\}$ then we obtain the classical Kaczmarz relaxation scheme for the solution of the linear system $Lu=f$, and the iteration scheme is equivalent to a Gauss-Seidel relaxation for the system $L^T L u = L^T f$.
In our case, one iteration of the alternative relaxation we want to study is composed as follows: we perform a Gauss-Seidel sweep in the interior of the domain, followed by $\lambda$ Kaczmarz iterations over ghost points and inner points close to the boundary (say within $\delta$ distance from the boundary). 

Another alternative is represented by the block relaxation~\cite{Dendy:Box}. As we point out in~\cite{CocoRusso:Elliptic}, the elimination of the boundary conditions is hard to perform in high dimensions, while in one dimension it is a trivial task and leads to a diagonally dominant linear system. A middle ground between the elimination of the boundary conditions and the relaxation operator we use in this paper is the block relaxation. Let us describe it in details. For each grid point $P \in \Omega_h \cup \Gamma_h$ we choose a stencil $St_P \subseteq \Omega_h \cup \Gamma_h$. For instance, if $P \in \Omega_h$ we choose $St_P=St_{P9} \cap \left( \Omega_h \cup \Gamma_h \right)$, where $St_{P9}$ is the $3 \times 3$ stencil centered at $P$, else if $P \in \Gamma_h$ we choose the stencil $St_P$ defined in \eqref{St_G}. One iteration of the alternative relaxation is composed as follows. We perform a Gauss-Seidel sweep in the interior of the domain except in grid points within $\delta$ distance from the boundary.
For each grid point $P \in \Omega_h \cup \Gamma_h$ within $\delta$ distance from the boundary we rewrite the linear system $Lu=f$ as follows (by a permutation of rows):
\[
 Lu=f \Leftrightarrow
\left(
\begin{array}{cc}
A_{1,1} & A_{1,2} \\
A_{2,1} & A_{2,2}
\end{array}
\right)
\left(
\begin{array}{c}
u_1 \\
u_2
\end{array}
\right)
=
\left(
\begin{array}{c}
f_1 \\
f_2
\end{array}
\right)
\]
where $u_1$ is referred to those grid points belonging to $St_P$. 
Therefore, we update the values of $u_1$ as:
\[
u_1 = A_{1,1}^{-1} \left( f_1-A_{1,2} u_2 \right).
\]

We perform a comparison between the relaxation proposed in this paper (that we call \textit{new iteration} in the following plots) and the two alternative relaxation described above. Such a comparison is carried out in terms of smoothing factor and convergence factor.
We perform the comparison using the TGCS for the test case of the circular domain \ref{NumTests:circle} with $N=64$.

\begin{figure}[!hbt]
 \begin{minipage}[c]{0.49\textwidth}
   	\centering
   	\captionsetup{width=0.70\textwidth}
		\includegraphics[width=0.90\textwidth]{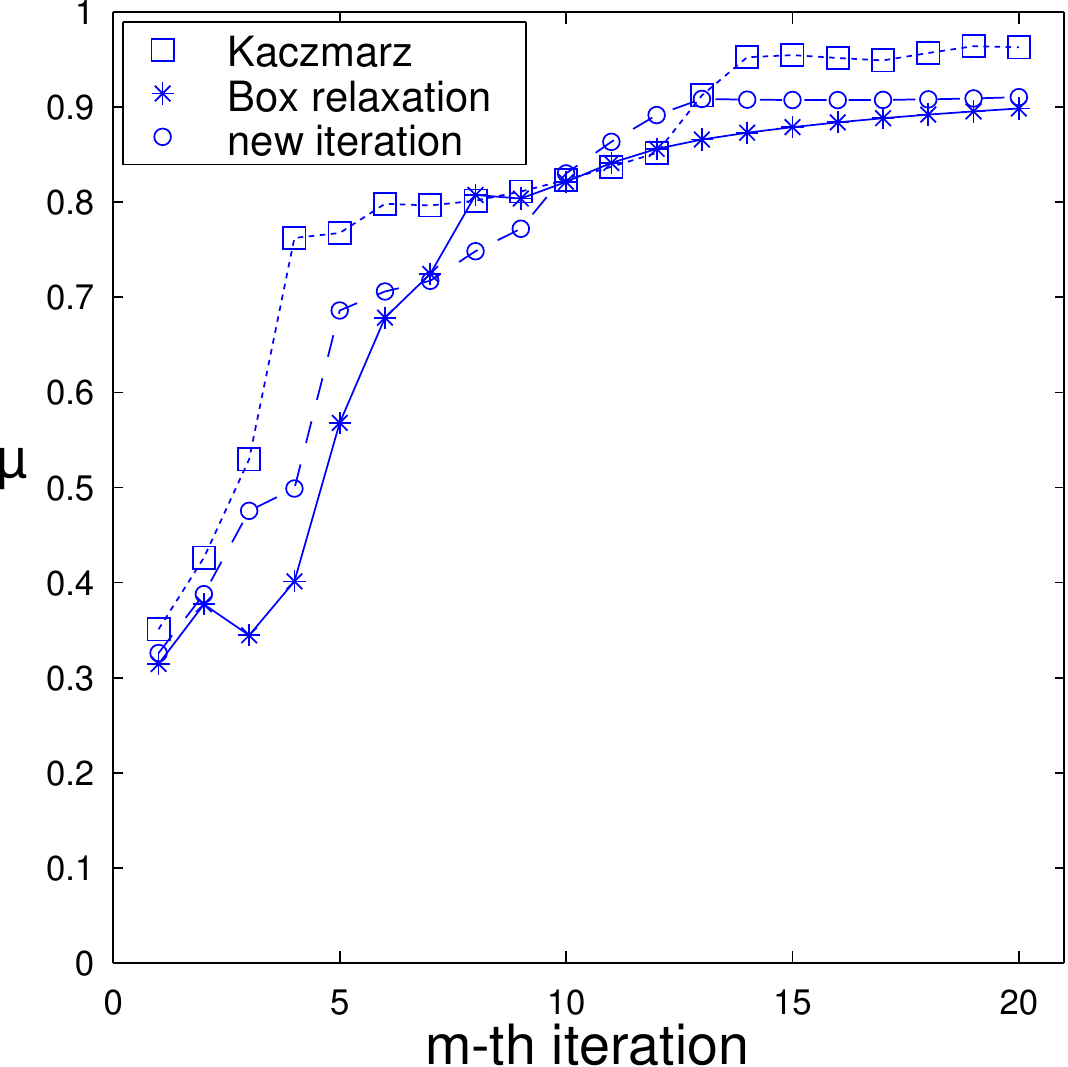}
	\caption{\footnotesize{Smoothing factor $\mu$ against the number of iterations for the three iterative schemes: Kaczmarz relaxation, Block relaxation, new iteration. }}
	\label{fig:smoothFact}
 \end{minipage}
 \begin{minipage}[c]{0.49\textwidth}
  	\centering
  	\captionsetup{width=0.80\textwidth}
		\includegraphics[width=0.90\textwidth]{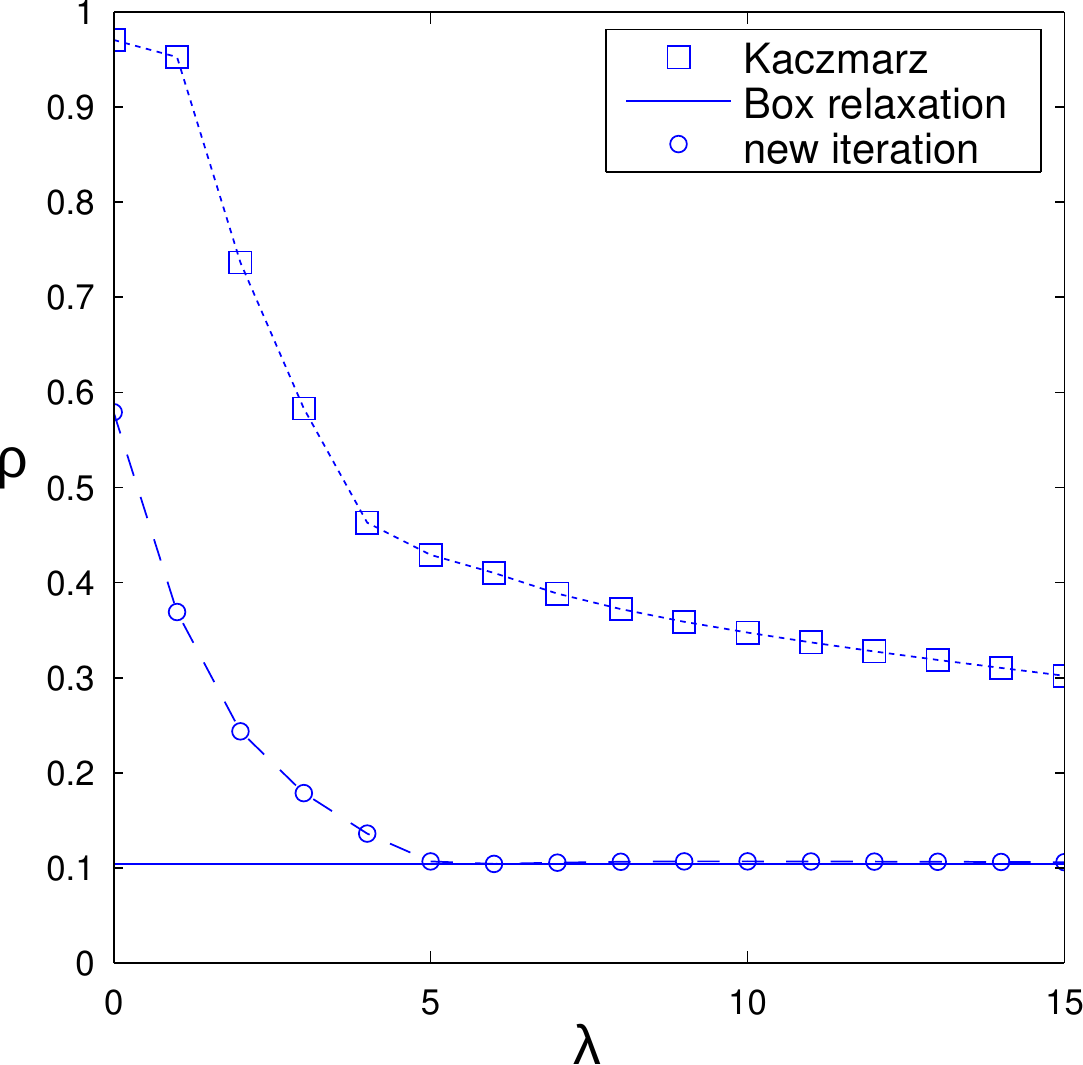}
	\caption{\footnotesize{Convergence factor $\rho$ of the entire multigrid against the number of extra-relaxations $\lambda$ for the three iterative schemes: Kaczmarz relaxation, Block relaxation, new iteration. }}
	\label{fig:convFactLambda}
 \end{minipage}
\end{figure}

In Fig.\ \ref{fig:smoothFact} we plot the smoothing factor $\mu$ for the three iteration schemes, which is estimated by the ratio of subsequent defects after each iteration, i.e.\
\[
\mu^{(m)}=\frac{\left\| \rh^{(m)} \right\|_{\infty}}{\left\| \rh^{(m-1)} \right\|_{\infty}}.
\]
In practice, we perform only the iteration schemes, without taking into account the effects of the multigrid procedure.
In order to better capture the behavior of the smoothing factor, we choose an initial guess being highly oscillant, for example $u=\sin(40\pi x) \sin(50 \pi y)$.

In Fig.\ \ref{fig:convFactLambda} we depict the convergence factor $\rho$ for the Kaczmarz and the new iteration against the number of extra-relaxations $\lambda$ (for comparison, we also plot the convergence factor of the block relaxation as an horizontal line, since it does not depend on $\lambda$). After five extra-relaxations, the new iteration reaches a plate configuration, since it achieves the convergence factor of the Gauss-Seidel smoother for inner equations, i.e.\ the convergence factor predicted by the LFA (see Table \ref{table:rholoc}). The Kaczmarz iteration falls down slower, while the block iteration already provides the optimal convergence factor.
The computational cost of five point-iterations of the new method is considerably lower than the cost of one block-iteration.

\subsubsection{Ellipsoidal domain}\label{NumTests:ellipse}
The level-set function is:
\[
\phi(x,y)=\frac{(X(x,y)-\sqrt{2}/20)^2}{0.563^2}+\frac{(Y(x,y)-\sqrt{3}/30)^2}{0.263^2}-1
\]
where
\[
 X(x,y)=\cos(\pi/6)\,x-\sin(\pi/6)\,y, \; \; \; X(x,y)=\sin(\pi/6)\,x+\cos(\pi/6)\,y,
\]
and the zero level-set is represented in Fig. \ref{fig:domains} (top-right).
The convergence factor obtained are listed in Table \ref{table:rhoE} (for $\nu=\nu_1+\nu_2=2$ and $\nu=3$). We observe as the convergence factor degrade choosing a coarsest grid too much coarse, but starting from a certain level of coarsest grid it is relatively close to the predicted convergence factor by LFA (Table \ref{table:rholoc}).

 \begin{table}[!hbt]
\captionsetup{width=0.79\textwidth}
\caption{ \footnotesize{Convergence factor of the numerical test of Section \ref{NumTests:ellipse}. We use $N \times N$ number of grid points in the finest grid; $N_c \times N_c$ number of grid points in the coarsest grid. Left table: $\nu=\nu_1+\nu_2=2$, right table $\nu=\nu_1+\nu_2=3$. }} % title of Table
\begin{minipage}{0.45\textwidth}
\centering      % used for centering table
\begin{tabular}{|| c c || c | c | c | c | c ||}  % centered columns (4 columns)
\hline\hline                        %inserts double horizontal lines
 & $N$ & 16 & 32 & 64 & 128 & 256  \\ [0.5ex] % inserts table
$N_c$ & &  &  &  &  & \\ 
\hline\hline                        %inserts double horizontal lines
8 & & 0.34 & 0.09 & 0.14 & 0.14 & 0.15\\ 
16 & &  & 0.65 & 0.45 & 0.19 & 0.15\\ 
32 & &  &  & 0.14 & 0.14 & 0.15\\ 
64 & &  &  &  & 0.15 & 0.15\\ 
128 & &  &  &  &  & 0.15\\ 
\hline\hline                        %inserts double horizontal lines
 \end{tabular} 
\end{minipage} 
\begin{minipage}{0.45\textwidth}
\centering      % used for centering table
\begin{tabular}{|| c c || c | c | c | c | c ||}  % centered columns (4 columns)
\hline\hline                        %inserts double horizontal lines
 & $N$ & 16 & 32 & 64 & 128 & 256  \\ [0.5ex] % inserts table
$N_c$ & &  &  &  &  & \\ 
\hline\hline                        %inserts double horizontal lines
8 & & 0.44 & 0.06 & 0.12 & 0.11 & 0.09\\ 
16 & &  & 0.55 & 0.30 & 0.09 & 0.09\\ 
32 & &  &  & 0.13 & 0.10 & 0.09\\ 
64 & &  &  &  & 0.12 & 0.08\\ 
128 & &  &  &  &  & 0.09\\ 
\hline\hline                        %inserts double horizontal lines
 \end{tabular} 
 \end{minipage}
  \label{table:rhoE}
 \end{table} 

\subsubsection{Saddle-shaped domain}\label{NumTests:saddle}
The level-set function is:
\[
\phi(x,y)=9\left(\frac{1}{2}x-\frac{\sqrt{3}}{2}y\right)^2+\left(\frac{3\sqrt{3}}{2}x+\frac{3}{2}y-1\right)^2\sin\left(\frac{3\sqrt{3}}{2}x+\frac{3}{2}y-1\right)-1
\]
and the zero level-set is represented in Fig. \ref{fig:domains} (bottom-left).
The convergence factor obtained for $\nu=\nu_1+\nu_2=3$ are listed in Table \ref{table:rhoSF3} (left). Also in this case, only if we choose $N=16$ and $N_c=8$ (which actually is TGCS) the convergence factor is degraded. 

 \begin{table}[!hbt]
\captionsetup{width=0.79\textwidth}
\caption{ \footnotesize{Convergence factor of the numerical test of Sections \ref{NumTests:saddle} (left) and \ref{NumTests:flower} (right). We use $N \times N$ number of grid points in the finest grid; $N_c \times N_c$ number of grid points in the coarsest grid. In this test we use $\nu=\nu_1+\nu_2=3$. For the flower-shaped domain, if $N_c=8$ the multigrid does not converge. }} % title of Table
\begin{minipage}{0.45\textwidth}
\centering      % used for centering table
\begin{tabular}{|| c c || c | c | c | c | c ||}  % centered columns (4 columns)
\hline\hline                        %inserts double horizontal lines
 & $N$ & 16 & 32 & 64 & 128 & 256  \\ [0.5ex] % inserts table
$N_c$ & &  &  &  &  & \\ 
\hline\hline                        %inserts double horizontal lines
8 & & 0.36 & 0.08 & 0.09 & 0.12 & 0.09\\ 
16 & &  & 0.12 & 0.09 & 0.12 & 0.09\\ 
32 & &  &  & 0.09 & 0.12 & 0.09\\ 
64 & &  &  &  & 0.13 & 0.09\\ 
128 & &  &  &  &  & 0.09\\ 
\hline\hline                        %inserts double horizontal lines
 \end{tabular} 
\end{minipage} 
\begin{minipage}{0.45\textwidth}
\centering      % used for centering table
\begin{tabular}{|| c c || c | c | c | c | c ||}  % centered columns (4 columns)
\hline\hline                        %inserts double horizontal lines
 & $N$ & 16 & 32 & 64 & 128 & 256  \\ [0.5ex] % inserts table
$N_c$ & &  &  &  & & \\ 
\hline\hline                        %inserts double horizontal lines
8 & & n.c. & n.c. & n.c. & n.c. & n.c.\\ 
16 & & & 0.89 & 0.75 & 0.50 & 0.25\\ 
32 & & &  & 0.49 & 0.25 & 0.12\\ 
64 & & &  &  & 0.24 & 0.11\\ 
128 & & &  &  &  & 0.09\\ 
\hline\hline                        %inserts double horizontal lines
 \end{tabular}  
 \end{minipage}
  \label{table:rhoSF3}
 \end{table} 

\subsubsection{Flower-shaped domain}\label{NumTests:flower}
The level-set function is:
\[
\phi = r-0.5-\frac{y^5+5 x^4y-10 x^2 y^3}{5 r^5}, \; \; \; r=\sqrt{x^2+y^2}
\]
and the zero level-set is represented in Fig. \ref{fig:domains} (bottom-right).
The convergence factor obtained for $\nu=\nu_1+\nu_2=3$ are listed in Table \ref{table:rhoSF3} (right). This is the hardest numerical test, because of the indentation of the boundary. We need to start from a coarsest level $N_c=32$ to correctly capture the boundary profile and to make the discretization accurate.
  
\section*{Conclusion}
A multigrid technique for Poisson equation on an arbitrary domain and mixed boundary conditions is presented. This multigrid strategy can be applied to a general framework of ghost-point method in a regular Cartesian grid, in case of non-eliminated boundary conditions. Suitable transfer operators for inside equations and boundary conditions are provided. The convergence rate is improved by adding some extra-relaxations on the ghost points and in a narrow band of inside grid points close to the boundary.
Numerical tests on different geometries have been performed. On simple domains (such as circle or ellipse) the optimal convergence factor is reached even if the method is used on very coarse grids, while for more complex domains (such as the flower-shaped one) the optimal convergence factor is obtained only on sufficiently fine grids.
A comparison with other treatments of the boundary condition smoothing procedure has been carried out (Kaczmarz and Block relaxation), confirming that the smoother proposed in this paper is better in terms of convergence factor, and not worse in terms of smoothing factor.

The application of this method to solve the pressure equation coming from the projection method of Chorin~\cite{Chorin:projection, Chorin:projection1997} in the framework of the incompressible Navier-Stokes equation is in preparation.
Several extensions of the discretization technique and multigrid approach are presently under investigation. We mention the case of discontinuous coefficients, which models, for example, a system composed by different materials separated by an interface; in such a case the method is suitably modified in order to achieve second order accuracy and a convergence factor being independent on the jump in the coefficient. A preliminary result can be found in~\cite{CocoRusso:Ischia2010}.
Another extension concerns the convection-diffusion equation in a moving domain, in order to study applications modeled by a Stefan-type problem.
All these extensions will be coupled with the use of Adaptive Mesh Refinement to obtain accurate solution in the case of domain with complex boundary.

%\tableofcontents

\addcontentsline{toc}{chapter}{References}
\bibliographystyle{abbrv}
\bibliography{bibliography}

\begin{thebibliography}{99}

\bibitem{Li:IIM_MG}
L.~Adams and Z.~Li.
\newblock {The immersed interface/multigrid methods for interface problems.}
\newblock {\em Journal of Scientific Computing}, 24:463--479, 2002.

\bibitem{Bramble:ell}
J.~H. Bramble and B.~E. Hubbard.
\newblock {Approximation of solutions of mixed boundary value problems for
  Poisson's equation by finite differences}.
\newblock {\em J. Assoc. Comput. Mach.}, 12:114--123, 1965.

\bibitem{Brandt:RigorousAnalysisMG}
A.~Brandt.
\newblock {Rigorous Quantitative Analysis of Multigrid, I: Constant
  Coefficients Two-Level Cycle with L2-Norm}.
\newblock {\em SIAM Journal on Numerical Analysis}, 31:1695--1730, 1994.

\bibitem{Briggs:MG}
W.~L. Briggs, V.~E. Henson, and S.~F. McCormick.
\newblock {\em {A Multigrid Tutorial}}.
\newblock SIAM, 2000.

\bibitem{Caflish:IslandDynamics}
R.~E. Caflisch, M.~F. Gyure, B.~Merriman, S.~J. Osher, C.~Ratsch, D.~D.
  Vvedensky, and J.~J. Zinck.
\newblock {Island dynamics and the level set method for epitaxial growth}.
\newblock {\em Applied Mathematics Letters}, 4:13--22, 1999.

\bibitem{Catalano:MG}
L.~A. Catalano, A.~Dadone, V.~S.~E. Daloiso, and D.~Scardigno.
\newblock {A multigrid procedure for Cartesian ghost-cell methods}.
\newblock {\em International Journal for Numerical Methods in Fluids},
  58:743--750, 2008.

\bibitem{Iollo:penalization}
F.~Chantalat, C.-H. Bruneau, C.~Galusinski, and A.~Iollo.
\newblock {Level-set, penalization and cartesian meshes: A paradigm for inverse
  problems and optimal design}.
\newblock {\em Journal of Computational Physics}, 228:6291--6315, 2009.

\bibitem{Gibou:quadtree}
H.~Chen, C.~Min, and F.~Gibou.
\newblock {A supra-convergent finite difference scheme for the Poisson and heat
  equations on irregular domains and non-graded adaptive Cartesian grids}.
\newblock {\em Journal of Scientific Computing}, 31:19--60, 2007.

\bibitem{vanLeer:CFLCutCell}
Y.~Chiang, B.~V. Leer, and K.~G. Powell.
\newblock {Simulation of unsteady inviscid flow on an adaptively refined
  Cartesian grid}.
\newblock In {\em AIAA Paper}, 1999.

\bibitem{Chorin:projection}
A.~Chorin.
\newblock {Numerical solution of the NavierÃ¢â¬âStokes Equations}.
\newblock {\em Mathematics of Computation}, 22:745--762, 1968.

\bibitem{Chorin:projection1997}
A.~Chorin.
\newblock {A numerical method for solving incompressible viscous flow
  problems}.
\newblock {\em Journal of Computational Physics}, 135:115--125, 1997.

\bibitem{Clarke:CFLCutCell}
D.~Clarke, M.~Salas, and H.~Hassan.
\newblock {Euler calculations for multielement airfoils using Cartesian grids}.
\newblock {\em AIAA Jounal}, 24:353--358, 1986.

\bibitem{CocoRusso:Elliptic}
A.~Coco and G.~Russo.
\newblock {A fictitious time method for the solution of Poisson equation in an
  arbitrary domain embedded in a square grid}.
\newblock {\em Journal of Computation Physics}.
\newblock Under revision.

\bibitem{CocoRusso:Ischia2010}
A.~Coco and G.~Russo.
\newblock {Second order multigrid methods for elliptic problems with
  discontinuous coefficients on an arbitrary interface, I: one dimensional
  problems}.
\newblock {\em Numerical Mathematics: Theory, Methods and Applications}.
\newblock Accepted.

\bibitem{CFL:CFL}
R.~Courant, K.~Friedrichs, and H.~Lewy.
\newblock {On the partial difference equations of mathematical physics}.
\newblock {\em IBM J. Res. Develop.}, 11:215--234, 1967.

\bibitem{Dendy:Box}
J.~Dendy.
\newblock {Black box multigrid.}
\newblock {\em Journal of Computational Physics}, 48:366--386, 1982.

\bibitem{Gibou:reinizialization}
A.~du~Ch\'en\'e, C.~Min, and F.~Gibou.
\newblock { Second-Order Accurate Computation of Curvatures in a Level Set
  Framework Using Novel High Order Reinitialization Schemes}.
\newblock {\em Journal of Scientific Computing archive}, 35:114--131, 2008.

\bibitem{Fedkiw:GFM}
R.~Fedkiw, T.~Aslam, B.~Merriman, and S.~Osher.
\newblock {A Non-Oscillatory Eulerian Approach to Interfaces in Multimaterial
  Flows (The Ghost Fluid Method)}.
\newblock {\em Journal of Computational Physics}, 152:457--492, 1999.

\bibitem{Gibou:Ghost}
F.~Gibou and R.~Fedkiw.
\newblock {A second-order-accurate symmetric discratization of the poisson
  equation on irregular domains}.
\newblock {\em Journal of Computational Physics}, 176:205--227, 2002.

\bibitem{Gibou:fourth_order}
F.~Gibou and R.~Fedkiw.
\newblock {A fourth order accurate discretization for the laplace and heat
  equations on arbitary domains, with applications to the stefan problem}.
\newblock {\em Journal of Computational Physics}, 202:577--601, 2005.

\bibitem{Hackbusch:MG}
W.~Hackbusch.
\newblock {\em {Multi-grid methods and applications}}.
\newblock Springer, 1985.

\bibitem{Hackbusch:elliptic}
W.~Hackbusch.
\newblock {\em {Elliptic Differential Equations: Theory and Numerical
  Treatment}}.
\newblock Springer, 2003.

\bibitem{Berger:rotated}
C.~Helzel, M.~J. Berger, and L.~R. J.
\newblock {A high-resolution rotated grid method for conservation laws with
  embedded geometries. (English summary)}.
\newblock {\em SIAM J. Sci. Comput.}, 26:785--809, 2005.

\bibitem{Colella:PoissonFV}
H.~Johansen and P.~Colella.
\newblock {A Cartesian Grid Embedded Boundary Method for Poisson Equation on
  Irregular Domains}.
\newblock {\em Journal of Computational Physics}, 147:60--85, 1998.

\bibitem{Sethian:level_set}
J.Sethian.
\newblock {\em { Level Set Methods and Fast Marching Methods: Evolving
  Interfaces in Computational Geometry, Fluid Mechanics, Computer Vision and
  Materials Science}}.
\newblock Cambridge University Press, 1999.

\bibitem{Kaczmarz:Kiteration}
S.~Kaczmarz.
\newblock {Angenäherte Auflösung von Systemen linearer Gleichungen.}
\newblock {\em Bulletin International de l'Académie Polonaise des Sciences et
  des Lettres. Classe des Sciences Mathématiques et Naturelles. Série A,
  Sciences Mathématiques}, 35:355--357, 1937.

\bibitem{LeVequeLi:IIM}
R.~LeVeque and Z.~Li.
\newblock {The immersed interface method for elliptic equations with
  discontinuous coefficients and singular sources}.
\newblock {\em SIAM J. Numer. Anal.}, 31:1019--1044, 1994.

\bibitem{Lonher:ImmersedEmbedded}
R.~L{\"o}hner, J.~R. Cebral, F.~E. Camelli, S.~Appanaboyina, J.~D. Baum, E.~L.
  Mestreau, and O.~A. Soto.
\newblock {Adaptive embedded and immersed unstructured grid techniques}.
\newblock {\em Comput. Methods Appl. Mech. Engrg.}, 197:217--2197, 2008.

\bibitem{Ma:MG}
Z.~H. Ma, L.~Qian, D.~M. Causon, H.~B. Gu, and C.~G. Mingham.
\newblock {A Cartesian ghost-cell multigrid Poisson solver for incompressible
  flows}.
\newblock {\em International Journal for Numerical Methods in Engineering},
  85:230--246, 2011.

\bibitem{Nochetto:StefanProblemFEM}
R.~H. Nochetto, M.~Paolini, and C.~Verdi.
\newblock {An adaptive finite element method for two-phase stefan problem in
  two space dimensions. part ii: Implementation and numerical experiments}.
\newblock {\em SIAM J. Sci. Stat. Comput.}, 12:1207--1244, 1991.

\bibitem{Osher-Fedkiw:level_set}
S.~Osher and R.~Fedkiw.
\newblock {\em {Level Set Methods and Dynamic Implicit Surfaces}}.
\newblock Springer-Verlag New York, Applied Mathematical Sciences, 2002.

\bibitem{Gibou:Robin}
J.~Papac, F.~Gibou, and C.~Ratsch.
\newblock {Efficient Symmetric Discretization for the Poisson, Heat and
  Stefan-Type Problems with Robin Boundary Conditions}.
\newblock {\em Journal of Computational Physics}, 229:875--889, 2010.

\bibitem{Peskin:immersedInterface}
C.~S. Peskin.
\newblock {Numerical analysis of blood flow in the heart}.
\newblock {\em Journal of Computational Physics}, 25:220--252, 1977.

\bibitem{Quart:PDE}
A.~Quarteroni.
\newblock {\em {Numerical models for differential problems}}.
\newblock Springer, 2009.

\bibitem{QuartSacco:PDE}
A.~Quarteroni and R.~Sacco.
\newblock {\em {Numerical approximation of partial differential equations}}.
\newblock Springer, 1994.

\bibitem{Russo-Smereka:reconstruction}
G.~Russo and P.~Smereka.
\newblock {A remark on computing distance functions}.
\newblock {\em Journal of Computational Physics}, 163:51--67, 2000.

\bibitem{Schmidt:dendritiesFEM}
A.~Schmidt.
\newblock {Computation of three dimensional dendrites with finite elements}.
\newblock {\em Journal of Computational Physics}, 125:293--312, 1996.

\bibitem{Shortley-Weller:discretization}
G.~H. Shortley and R.~Weller.
\newblock {The numerical solution of laplace's equation.}
\newblock {\em J. Appl. Phys.}, 9:334--348, 1938.

\bibitem{Strikwerda:FD}
J.~C. Strikwerda.
\newblock {\em {Finite Difference Schemes and Partial Difference Equations.
  Second Edition.}}
\newblock SIAM, 2004.

\bibitem{SSO:level_set}
M.~Sussman, P.~Smereka, and S.~Osher.
\newblock {A level set approach for computing solutions to incompressible
  2-phase flow}.
\newblock {\em Journal of Computational Physics}, 114:146--159, 1994.

\bibitem{Trottemberg:MG}
U.Trottemberg, C.~Oosterlee, and A.~Schuller.
\newblock {\em {Multigrid}}.
\newblock Academic Press, 2000.

\bibitem{LeVequeLiWiegmann:CRACK}
A.~Wiegmann, Z.~Li, and R.~LeVeque.
\newblock {Crack Jump Conditions for Elliptic Problems}.
\newblock {\em Applied Mathematics Letters}, 12:81--88, 1999.

\end{thebibliography}
%\nocite{*}

\end{document}